\documentclass[10pt, final]{article}
\usepackage{color}
\usepackage{tikz}
\usetikzlibrary{arrows}

\usepackage{latexsym}
\usepackage[utf8]{inputenc}
\usepackage{srcltx}
\usepackage{amscd}
\usepackage{graphicx}
\usepackage{epsfig}
\usepackage{subfig}
\usepackage{amsmath}
\usepackage{amssymb}
\usepackage{mathrsfs}
\input xy
\xyoption{all}
\usepackage{srcltx}
\usepackage{float}
\usepackage{booktabs}
\usepackage{amssymb}
\usepackage{fixltx2e}
\usepackage{graphicx}

\usepackage{geometry}

\newcommand{\Z}{\mathbb{Z}}                     
\newcommand{\R}{\mathbb{R}}                     

\newcommand{\finedim}{\hfill $\Box$\\}            

\newcommand{\proof}{{\sl Proof.}\hspace{5pt}}   

\newtheorem{thm}{\sc Theorem}[section]      
\newtheorem{cor}[thm]{\sc Corollary}        
\newtheorem{lem}[thm]{\sc Lemma}            
\newtheorem{prop}[thm]{\sc  Proposition}     
\newtheorem{defn}[thm]{\sc Definition}      
\newtheorem{rem}[thm]{\sc Remark}       
\newtheorem{example}[thm]{\sc Example}           

\title{Modes in modern music from a topological viewpoint}
\author{Mattia G. Bergomi, Alessandro Portaluri}

\begin{document}
\date{\today}
\maketitle

\footnotetext[2]{A. Portaluri was supported by the grant PRIN2009 ``Critical
Point Theory and Perturbative Methods
for Nonlinear Differential Equations''.}

\begin{abstract}

The aim of this paper is twofold: on one side we review the classical concept
of 
{\em musical mode\/} from the viewpoint of modern music, reading it as a superimposition of a 
{\em base-chord\/} (seventh chord) and a {\em tension-chord\/} (triad). We  
associate to each modal scale an oriented plane graph whose homotopy properties give a
measure of 
the complexity of the base-chord associated to a certain mode. Using these graphs we prove the existence of
{\em special 
modes\/} which are not deducible in the standard way. 

On the other side we give a more deep musical insight by developing a 
braid theoretical interpretation of some cadential harmonic progressions in modern music  
and we use  braid theory in order to represent them and voice leadings among them.

A striking application is provided by the analysis of an harmonic fragment from  
{\em Peru\/} by {\em Tribal Tech\/}. 
We approximate the octatonic scale used in the improvisation by Scott 
Henderson, through the special mode {\em  $myxolydian\,\flat 2 \sharp 4$\/} and we finally associate a braid 
representation to the fragment we analysed. 

\end{abstract}

\maketitle

\tableofcontents

\section*{Introduction}

Classical harmony has been largely  studied, however musicians
often are hard put to understand the real nature of modes and the
concept of sonority in modern music. In Western music theory, the word ``mode'' (from Latin modus, 
`` measure, way, size, method'')
generally refers to a type of scale,  coupled with a set of 
characteristic melodic behaviours. More precisely
the concept of mode incorporates the idea of the diatonic scale, but
differs from 
it involving an element of melody type. This concerns particular  
repertories of short musical figures or groups of tones within a certain scale
so that, 
depending on the point of view, mode takes the meaning of either a 
``particularized scale'' or a ``generalized tune''. 

The goal of this paper is to construct modes from a genuine harmonic viewpoint; 
more precisely, a mode is not treated only like a scale (and hence an ordered sequence of
pitches) but 
also as a superimposition of two chords: a {\em base-chord\/} and a {\em tension
chord\/}. 
From a technical musical viewpoint, the base-chord is a seventh chord and the tension one is a triad. Of course as one
can easily 
figure out not every seventh chord is admissible as base-chord as well as not
every triad is 
compatible with a fixed seventh chord. Specifically, the base-chords have to be
into 7 classes corresponding to the types of chord we can deduce harmonizing each degree of the major, harmonic 
and melodic minor scales. Each base-chord and the fact that we deduce modes using three scales ``force'' the superimposition of some
specific triads which are chosen decomposing the classical modal scales; 
namely, not every triad can be used to produce a mode, as we shall see later. 

Starting from a fixed base-chord we shall consider only the triads that give us 
the classical modes and we shall associate in a quite elementary  way an
oriented 
planar connected graph $G$ (the order being established by the scale's degrees)
and 
hence a topological measure of complexity given by the $1- \chi(G)$ where
$\chi(\cdot)$ 
denotes the Euler characteristic. It is a well-known fact that the fundamental
group of $G$, 
i.e. $\pi_1(G)$ is a free group over $1-\chi(G)$ generators. Moreover
$\chi(\cdot)$ is  a
topological invariant.  Each mode is a subgraph of the associated graph 
such that each 
vertex (except the first and the last one) has valence at least $2$ (i.e. at
least two edges connect 
the vertex in the graph). 

Our first result is that out of the classical $21$ modes which we shall refer as
{\em 
standard modes\/}, other $12$  ``exceptional'' modes arise. 
We shall refer to them as {\em admissible special modes\/} (See proposition \ref{prop:specmodes}). The most important fact is that one of these
modes which is called $Myxolidian\,\flat 2\sharp 4,$
allows us to well ``approximate'' an octatonic half step/whole step scale used in a famous fusion piece, with a 7 notes scale, 
conserving entirely the notes of the chord played in the harmonic structure (a dominant chord in the case we examine). This fact
allows either the composer or the improviser to constantly see the notes of resolutions which are the base-chord and to clearly 
distinguish them from the tension notes which are collected in the tension-triad. Moreover in the $Myxolidian\,\flat 2\sharp 4$ we have
the presence of the tritone which makes the dominant chord recognizable, as the $\flat 2$ and the augmented fourth  which greatly 
contribute to the tension of the
mode.

We would like to stress the fact that modes attracted the attention of many
music theorists over 
the years. However, to the knowledge of the authors, the analysis carried over
is mainly of quantitative type and in fact it is based on
algebra 
and combinatorics. We shall refer to the papers of \cite{Agm89,Nol09,Row00, Maz02} and
references therein. 

Our approach is in the opposite direction being somehow qualitative in the
essence 
and for this reason the suitable mathematical tools are mainly based on the use 
of (elementary) algebraic topology. 

In contrast to the first part which is ``static'' in the essence, in the
second part of the paper we 
analyse some harmonic progressions. More precisely we take into account
some well-known cadential chord progressions in modern music and in the special case in jazz music which 
are based on perfect, deceptive, plagal and secondary dominant progressions. Our main
mathematical device in 
order to carry over our investigation is based on the geometrical interpretation
through braids which 
give us a suggestive and net interpretation of the mathematical structure behind
this corner of the musical 
world. We analyse the same harmonic fragment we used to test the special modes, which put on
evidence a sort of ``musical complexity'' through its representation with braids. 

This is the starting point of a forthcoming paper in which the 
braids structure of some cadences are used in order to construct suitable
topological spaces having 
highly non-trivial homotopic properties. These properties allow us to define in a
very clear and rigorous 
way a sort of {\em musical complexity\/} which could be used in order to
classify some musical pieces, as well as to distinguish the degree of complexity of a 
composition or improvisation. (See section \ref{sec:conclusions} for further details and future projects).


\section{Standard modes and superimposition of chords}\label{sec:grafi}

The aim of this section is to recall some well-known facts and definitions about
the classical 
construction of modes in modern music. A definition of mode \cite{Lev95} is that
it is a seven-note scale created by 
starting on any of the seven note of a major or a melodic minor scale.
The section is structured as follows: in the first paragraph we deduce the modes
 from the major, melodic minor and 
harmonic minor scale. In the second paragraph a  mathematical way to represent
chords and modes is shown. We conclude the 
section showing a non-standard way to think about modes as a superimposition of
two chords.

\subsection{Deducing the standard modes}

Following the definition given in \cite{Lev95}, we assume that a mode is an
heptatonic scale, however we include the harmonic 
minor in the set of scales we use to deduce the standard modes. The reason which
lead us to consider this scale is that it allows 
us to define naturally some really used and interesting modes like the
\emph{phrygian dominant}.
The classical way to deduce modes from scales is to consider each degree of the
scale as the root of the modal scale. 
In table \ref{tab:clasmode} we list the $21$ modes given by the major, melodic
minor and harmonic minor scale, examples 
have been built on the C major ($C,D,E,F,G,A,B$), melodic minor ($C,D,E\flat
,F,G,A,B$) and harmonic
minor ($C,D,E\flat ,F,G,A\flat ,B$) scale respectively.
\begin{table}[H]
\begin{centering}
\begin{tabular}{|c|c|c|c|}
\hline 
Scale & Degree & Modes & Example \\
\hline 
\hline 
Major & $I$ & Ionian &  {\small C - D - E - F - G - A - B}\\
& $II$ & Dorian & {\small D - E - F - G - A - B - C}\\
& $III$ & Phrygian  & {\small E - F - G - A - B - C - D}\\
& $IV$ & Lydian & {\small F - G - A - B - C - D - E}\\
& $V$ & Mixolydian & {\small G - A - B - C - D - E - F}\\
& $VI$ & Eolian & {\small A - B - C -D - E - F - G}\\
& $VII$ & Locrian & {\small  B - C - D - E - F - G - A}\\
\hline
\hline
Melodic Minor & $I$ & Hypoionian & {\small C - D - E$\flat$ - F - G - A - B}\\
& $II$ & Dorian $\flat 2$ & {\small D - E$\flat$ - F - G - A - B - C}\\
& $III$ & Lydian Augmented & {\small  E$\flat$ - F - G - A - B - C - D}\\
& $IV$ & Lydian Dominant & {\small  F - G - A - B - C - D - E$\flat$}\\
& $V$ & Mixolydian $\flat 13$ & {\small G - A - B - C - D - E$\flat$ - F}\\
& $VI$ & Locrian $\sharp 2$ & {\small A - B - C -D - E$\flat$ - F - G}\\
& $VII$ & Super locrian & {\small B - C - D - E$\flat$ - F - G - A}\\
\hline
\hline
Harmonic Minor & $I$ & Hypoionian $\flat 6$ & {\small C - D - E$\flat$ - F - G -
A$\flat$ - B}\\
& $II$ & Locrian $\sharp 6$ & {\small D - E$\flat$ - F - G - A$\flat$ - B - C}\\
& $III$ & Ionian augmented & {\small  E$\flat$ - F - G - A$\flat$ - B - C - D}\\
& $IV$ & Dorian $\sharp 4$ &{\small F - G - A$\flat$ - B - C - D - E$\flat$}\\
& $V$ & Phrygian dominant & {\small G - A$\flat$ - B - C - D - E$\flat$ - F}\\
& $VI$ & Lydian $\sharp 2$ & {\small A$\flat$ - B - C - D - E$\flat$ - F - G}\\
& $VII$ & Ultra locrian & {\small B - C - D - E$\flat$ - F - G - A$\flat$}\\
\hline 
\end{tabular}
\caption{The $21$ modes derived from the major, melodic minor and harmonic minor
scale.}
\label{tab:clasmode}
\par\end{centering}
\end{table}
Musically speaking, from the classification given in table \ref{tab:clasmode} it
is really hard to understand 
what a mode is: a non-trained listener could say that the dorian mode is simply
a $C$ major scale played from the 
second degree, but every musician knows that this is a great reduction. However
suffice to do a simple experiment 
to understand that modes are really distinguishable when they are played on the
seventh chord built on the degree 
of the scale associated to the modal scale. (See table \ref{tab:harmonization}). For a
non-trained listener a $D$-dorian
scale played on a $Cmaj7$ chord is simply a $C$ major scale played on its second
degree, but the same scale played 
on a $D-7$ chord sounds is perceived as a new scale, that results completely different
from the $C$-major scale.\\
Getting to the point, we focus our attention on this strong relation between the
modal scale and the chord we associate
to the scale, see table \ref{tab:clasmode} and \ref{tab:harmonization}.
\begin{table}[H]
\begin{centering}
\begin{tabular}{|c|c|c|c|}
\hline 
Scale & Degree & $7^{th}$ Chord & Arpeggio (example) \\
\hline 
\hline 
Major & $I$ & $maj7$ &  {\small C - E - G - B}\\
& $II$ & $-7$ & {\small D - F - A - C}\\
& $III$ & $-7$  & {\small E - G - B - D}\\
& $IV$ & $maj7$ & {\small F - A - C - E}\\
& $V$ & $7$ & {\small G - B - D - F}\\
& $VI$ & $-7$ & {\small A - C - E - G}\\
& $VII$ & $-7^{\flat 5}$ & {\small  B - D - F - A}\\
\hline
\hline
Melodic Minor & $I$ & $-maj7$ & {\small C - E$\flat$ - G - B}\\
& $II$ & $-7$ & {\small D - F - A - C}\\
& $III$ & $maj7^{\sharp 5}$ & {\small  E$\flat$ - G - B - D}\\
& $IV$ & $7$ & {\small  F - A - C - E$\flat$}\\
& $V$ & $7$ & {\small G - B - D - F}\\
& $VI$ & $-7^{\flat 5}$ & {\small A - C - E$\flat$ - G}\\
& $VII$ & $-7^{\flat 5}$ & {\small B - D - F - A}\\
\hline
\hline
Harmonic Minor & $I$ & $-maj7$ & {\small C - E$\flat$ - G - B}\\
& $II$ & $-7^{\flat 5}$ & {\small D - F - A$\flat$ - C}\\
& $III$ & $maj7^{\sharp 5}$ & {\small  E$\flat$ - G - B - D}\\
& $IV$ &  $-7$ &{\small F - A$\flat$ - C - E$\flat$}\\
& $V$ & $7$ & {\small G - B - D - F}\\
& $VI$ & $maj7$ & {\small A$\flat$ - C - E$\flat$ - G}\\
& $VII$ & $^{\circ7}$ & {\small B - D - F - A$\flat$}\\
\hline 
\end{tabular}
\caption{Seventh chord harmonization on the major, melodic minor and harmonic
minor scale. Chords notation: $Rmaj7$ 
is given by the root $R$, its major third ($M3$), perfect fifth ($P5$) and major
seventh ($M7$); we denote with $R-7$ 
the chord given by $R-m3-P5-m7$ where $m3$ and $m7$ are respectively the minor
third and the minor seventh in respect 
to the root note. $R-maj7$ is $R-m3-P5-m7$, $R-7^{\flat 5}$ is $R-m3-dim5-m7$,
$Rmaj7^{\sharp 5}$ is given by 
$R-M3-aug5-M7$ and $R^{\circ 7}$ is $R-m3-dim5-dim7$.}
\label{tab:harmonization}
\par\end{centering}
\end{table}


\subsection{Modes as a superimposition of chords}\label{par:modessuper}

In this paragraph we want to stress the importance of the harmonic choice which
lies behind the modal scale. 
Since either the chord and the modal scale is associated to a certain degree of
the scale, the notes of the chord have
to be part of the modal scale, more precisely they are the first, third, fifth
and seventh degree of the modal scale.

\begin{example}\label{example:base}
Consider $F$ lydian as a modal scale. The chord associated to this mode is
$Fmaj7$, then we have 
\begin{table}[H]
\begin{centering}
\begin{tabular}{c c}
$F$ lydian scale & $F-G-A-B-C-D-E$ \\
$Fmaj7$ arpeggio & $F-A-C-E$ \\
\end{tabular}
\par\end{centering}
\end{table} 
\end{example}
Every mode is associated to the proper seventh chord built on the root of the modal scale, 
this chord is called for brevity the \emph{base-chord}. We list them in table
\ref{tab:modchord}.
\begin{table}[H]
\begin{centering}
\begin{tabular}{|c|c|c|}
\hline 
Mode & Example - Base-chord & Example - Scale \\
\hline 
\hline 
Ionian & $Cmaj7$ & {\small C - D - E - F - G - A - B}\\
Dorian & $D-7$ & {\small D - E - F - G - A - B - C}\\
Phrygian  & $E-7$ & {\small E - F - G - A - B - C - D}\\
Lydian & $Fmaj7$ & {\small F - G - A - B - C - D - E}\\
Mixolydian & $G7$ & {\small G - A - B - C - D - E - F}\\
Eolian & $A-7$ & {\small A - B - C - D - E - F - G}\\
Locrian & $B-7^{\flat 5}$ & {\small  B - C - D - E - F - G - A}\\
\hline
\hline
Hypoionian & $C-maj7$ & {\small C - D - E$\flat$ - F - G - A - B}\\
Dorian $\flat 2$& $D-7$ & {\small D - E$\flat$ - F - G - A - B - C}\\
Lydian Augmented & $E\flat maj7^{\sharp 5}$ & {\small  E$\flat$ - F - G - A - B
- C - D}\\
Lydian Dominant & $F7$ & {\small  F - G - A - B - C - D - E$\flat$}\\
Mixolydian $\flat 13$ & $G7$ & {\small G - A - B - C - D - E$\flat$ - F}\\
Locrian $\sharp 2$ & $A-7^{\flat 5}$ & {\small A - B - C -D - E$\flat$ - F -
G}\\
Super locrian & $B-7^{b5}$ & {\small B - C - D - E$\flat$ - F - G - A}\\
\hline
\hline
Hypoionian $\flat 6$ & $C-maj7$ & {\small C - D - E$\flat$ - F - G - A$\flat$ -
B}\\
Locrian $\sharp 6$ & $D-7^{\flat 5}$ & {\small D - E$\flat$ - F - G - A$\flat$ -
B - C}\\
Ionian augmented & $E\flat maj7^{\sharp 5}$ & {\small  E$\flat$ - F - G -
A$\flat$ - B - C - D}\\
Dorian $\sharp 4$& $C-7$ & {\small F - G - A$\flat$ - B - C - D - E$\flat$}\\
Phrygian dominant & $G7$ & {\small G - A$\flat$ - B - C - D - E$\flat$ - F}\\
Lydian $\sharp 2$ & $Amaj7$ & {\small A$\flat$ - B - C -D - E$\flat$ - F - G}\\
Ultra locrian & $B^{\circ 7}$ & {\small B - C - D - E$\flat$ - F - G -
A$\flat$}\\
\hline 
\end{tabular}
\caption{The twenty-one modes derived from the major, melodic minor and harmonic
minor scale and their base-chords.}
\label{tab:modchord}
\par\end{centering}
\end{table}

From either example \ref{example:base} and table \ref{tab:modchord} one can also
deduce that deleting the 
base-chord from a modal scale a triad remains, we call that triad the
\emph{tension-triad} which is 
composed by the second, the fourth and the sixth degree of the modal scale.

\begin{example}
Considering the same setting of example \ref{example:base}, we have

\begin{table}[H]
\begin{centering}
\begin{tabular}{c c}
$F$ lydian scale & $F-G-A-B-C-D-E$ \\
Base-chord arpeggio & $F-A-C-E$ \\
Tension-traid arpeggio & $G-B-D$ \\
\end{tabular}
\par\end{centering}
\end{table} 
Thus it is possible to think about the $F$ lydian mode as the superimposition of
a $Fmaj7$ chord and a $G$ major triad.
\end{example}

Every modal scale can be decomposed uniquely in a seventh chord built on its
root note and a triad built on 
its second degree. As we said before the modal scale is recognizable if it is
played on its base-chord, so 
one can consider the base-chord as the set of stable notes and the tension-triad
as the collection of unstable 
notes for a certain mode, see table \ref{tab:superimp} for a complete
description of modes in terms of base-chords and tension-triads.

\begin{table}[H]
\begin{centering}
\begin{tabular}{|c|c|c|c|}
\hline 
Mode & Scale & Base-chord & Tension-triad \\
\hline 
\hline 
$C$ Ionian & {\small C - D - E - F - G - A - B} & $Cmaj7$ & $D-$\\
$D$ Dorian & {\small D - E - F - G - A - B - C} & $D-7$ & $E-$ \\
$E$ Phrygian & {\small E - F - G - A - B - C - D}  & $E-7$ & $F$ \\
$F$ Lydian  & {\small F - G - A - B - C - D - E} & $Fmaj7$ & $G$ \\
$G$ Mixolydian  & {\small G - A - B - C - D - E - F} & $G7$ & $A-$ \\
$A$ Eolian  & {\small A - B - C - D - E - F - G} & $A-7$ & $B-\flat 5$\\
$B$ Locrian  & {\small  B - C - D - E - F - G - A} & $B-7^{\flat 5}$ & $C$\\
\hline
\hline
$C$ Hypoionian  & {\small C - D - E$\flat$ - F - G - A - B} & $C-maj7$ & $D-$\\
$D$ Dorian $\flat 2$ & {\small D - E$\flat$ - F - G - A - B - C } & $D-7$ &
$E\flat^{\sharp 5}$\\
$E\flat$ Lydian Augmented  & {\small  E$\flat$ - F - G - A - B - C - D} &
$E\flat maj7^{\sharp 5}$ & $F$\\
$F$ Lydian Dominant & {\small  F - G - A - B - C - D - E$\flat$} & $F7$ & $G$\\
$G$ Mixolydian $\flat 13$ & {\small G - A - B - C - D - E$\flat$ - F} & $G7$ &
$A-\flat 5$\\
$A$ Locrian $\sharp 2$ & {\small A - B - C -D - E$\flat$ - F - G} & $A-7^{\flat
5}$ & $B -\flat 5$\\
$B$ Super locrian & {\small B - C - D - E$\flat$ - F - G - A} & $B-7^{b5}$ &
$C-$\\
\hline
\hline
$C$ Hypoionian $\flat 6$  & {\small C - D - E$\flat$ - F - G - A$\flat$ - B} &
$C-maj7$ & $D-\flat 5$\\
$D$ Locrian $\sharp 6$ & {\small D - E$\flat$ - F - G - A$\flat$ - B - C} &
$D-7^{\flat 5}$ & $E\flat^{\sharp 5}$\\
$E\flat$ Ionian augmented & {\small  E$\flat$ - F - G - A$\flat$ - B - C - D} &
$E\flat maj7^{\sharp 5}$ & $F$\\
$F$ Dorian $\sharp 4$ & {\small F - G - A$\flat$ - B - C - D - E$\flat$} & $C-7$
& $G$\\
$G$ Phrygian dominant & {\small G - A$\flat$ - B - C - D - E$\flat$ - F} & $G7$
& $A-\flat 5$\\
$A\flat$ Lydian $\sharp 2$ & {\small A$\flat$ - B - C -D - E$\flat$ - F - G} &
$Amaj7$ & $B-\flat 5$\\
$B$ Ultra locrian & {\small B - C - D - E$\flat$ - F - G - A$\flat$} & $B^{\circ
7}$ & $C-$\\
\hline 
\end{tabular}
\caption{Modes as a superimposition of two chords. Notation $^{\sharp 5}$ means
that the triad is 
augmented ($R$ - $M3$ - $aug 5$), while $^{\flat 5}$ means that it is diminished
($R$ - $m3$ - $dim5$).}
\label{tab:superimp}
\par\end{centering}
\end{table}

This decomposition which arise naturally from the standard way of deducing modes
from the classical 
scales, is the reason why we want to introduce in the next paragraph the space
of $4$ and $3$ note-chords.

\subsection{The space of chords: some notations}

The aim of this section is to recall some well-known facts about the space of 
chords and to fix our notations. Our basic references are \cite{CQT08,Tym11, Tym06, Tym09}  and 
references therein.

From a mathematical point of view, we recall that {\em pitches\/} are modeled by
real numbers. A {\em pitch class\/} (modulus octave) is modeled by a point in the 
$1$-dimensional torus $\mathbb T^1:=\R/(12
\Z)$. In general the space  $\mathbb T^n:= \R^n/(12 \Z)^n$ represents the space of ordered 
pairs of pitch
classes. (Cfr. \cite[Section 2-5]{Tym06}, for further details). However this space is bigger
than the space of $n$-notes 
chords which corresponds to the {\em unordered \/} pairs of pitch class. Thus by
passing to the quotient of 
$\mathbb T^n$ by the group of all permutations over the set of $n$ elements,
namely $\mathfrak S(n)$,
we get the space of $\mathfrak C_n$ of $n$-notes chords. More precisely:
\[
 \mathfrak C_n:=\mathbb T^n/\mathfrak S(n)= \R^n/(12\Z )^n\rtimes \mathfrak S(n)
\]
where $\rtimes$ denotes the semi-direct product. We observe that $\mathfrak C_2$
is the  M\"obius strip 
which is a not orientable surface. Also for dimension higher than two the quotient is
an orbifold. In the 
following we mainly consider the space of seventh chords (i.e. $4$-notes chords)
and triads ($3$-notes chords) that 
corresponds respectively to the orbifold constructed gluing with a twist two opposite top and bottom
faces of a $4$-dimensional 
and $3$-dimensional prism. (Cfr. \cite[Fig. S5-S6]{Tym06}). However we do not use the geometric structure of the space of chords except as an ambient space.

\begin{defn}\label{def:intersezione}
Given two chords $[A] \in\mathfrak C_n$ and $[C] \in\mathfrak C_m$, we define
the 
{\em chord intersection\/} of both as the chord $[D] \in\mathfrak C_k$ for $k
\leq \min\{m,n\}$ defined
by the maximal dimensional subset of notes appearing both in $[A]$ and $[C]$.
\end{defn}

These structures allow us to represent a {\em mode\/} as the superimposition of
a $4$-notes chord and a $3$-notes chord. 
More precisely:

\begin{defn}\label{def:modo}
We set:
\[
\mathscr M:= \mathfrak C_4 \times \mathfrak C_3. 
\]
A {\em musical mode\/} or simply a {\em mode\/} is a point $m=([B],[T])\in
\mathscr M$ 
such that $[B]\cap[T]=\emptyset$.
\end{defn}

As we will see later on, not every point in $\mathscr{M}$ is a mode.
The decomposition has be formalized, thus we can give another classification of
modal scales fixing the
base-chord $B$ and varying the tension-triad to associate every possible modal
choice to a certain seventh 
chord (available tension-triads are the ones we derived in paragraph
\ref{par:modessuper}, table \ref{tab:superimp}). 
See table \ref{tab:modclass} for a list of all possible modal choices on a fixed
type of seventh chord \cite{BG12}.

\begin{table}[H]
\begin{centering}
\begin{tabular}{|c|c|c|}
\hline 
Base-chord $maj7$ & Modes & Example (root C)\tabularnewline
\hline 
\hline 
{\small T-M III-P V-M VII} & Ionian & {\small C - D - E - F - G - A - B}\\
\cline{2-3} 
 & Lydian & {\small C - D - E - F$\sharp$ - G - A - B}\\
\cline{2-3} 
 & Lydian $\sharp2$ & {\small C - D$\sharp$ - E - F$\sharp$ - G - A - B}\\
\hline 
\hline 
Base-chord \emph{ $maj7^{\,\sharp5}$} & Modes & Example (root C)\\
\hline 
\hline 
{\small T-M III-aug V-M VII} & Lydian augmented & {\small C - D - E - F$\sharp$
- G$\sharp$ - A - B}\\
\cline{2-3} 
 & Ionian augmented & {\small C - D - E - F - G$\sharp$ - A - B}\\
\hline 
\hline 
Base-chord $7$ & Modes & Example (root C)\\
\hline 
\hline 
{\small T-III M-V P-VII m} & Mixolydian & {\small C - D - E - F - G - A -
B$\flat$}\\
\cline{2-3} 
 & Mixolydian $\flat13$ & {\small C - D - E - F - G - A$\flat$ - B$\flat$}\\
\cline{2-3} 
 & Phrygian dominant & {\small C - D$\flat$ - E - F - G - A$\flat$ - B$\flat$}\\
\cline{2-3} 
 & Lydian dominant & {\small C - D - E - F$\sharp$ - G - A - B$\flat$}\\
\hline 
\hline 
Base chord $-7$ & Modes & Example (root C)\\
\hline 
\hline 
{\small T-m III-P V-m VII} & Dorian & {\small C - D - E$\flat$ - F - G - A -
B$\flat$}\\
\cline{2-3} 
 & Phrygian & {\small C - D$\flat$ - E - F - G - A$\flat$ - B$\flat$}\\
\cline{2-3} 
 & Eolian & {\small C - D - E$\flat$ - F - G - A$\flat$ - B$\flat$}\\
\cline{2-3} 
 & Dorian $\flat2$ & {\small C - D$\flat$ - E$\flat$ - F - G - A - B$\flat$}\\
\cline{2-3} 
 & Dorian $\sharp4$\textcolor{white}{..} & {\small C - D - E$\flat$ - F$\sharp$
- G - A - B$\flat$}\\
\hline 
\hline 
Base chord $-7^{\,\flat5}$ & Modes & Example (root C)\\
\hline 
\hline 
{\small T-m III-dim V-m VII} & Locrian & {\small C - D$\flat$ - E$\flat$ - F -
G$\flat$ - A$\flat$ - B$\flat$}\\
\cline{2-3} 
 & Locrian $\sharp2$ & {\small C - D - E$\flat$ - F - G$\flat$ - A$\flat$ -
B$\flat$}\\
\cline{2-3} 
 & Superlocrian & {\small C - D$\flat$ - E$\flat$ - F$flat$ - G$\flat$ -
A$\flat$ - B$\flat$}\\
\cline{2-3} 
 & Locrian $\sharp6$ & {\small C - D$\flat$ - E$\flat$ - F - G$\flat$ - A -
B$\flat$}\\
\hline 
\hline 
Base-chord $-maj7$ & Modes & Example (root C)\\
\hline 
\hline 
{\small T-m III-P V-M VII} & Hypoionian & {\small C - D - E$\flat$ - F - G - A -
B}\\
\cline{2-3} 
 & Hypoionian $\flat6$ & {\small C - D - E$\flat$ - F - G - A$\flat$ - B}\\
\hline 
\hline 
Base-chord $\circ7$ & Modes & Example (root C)\\
\hline 
\hline 
{\small T-m III-dim V-dim VII} & Ultralocrian & {\small C - D$\flat$ - E$\flat$
- F - G$\flat$ - A$\flat$ - B$\flat\flat$}\\
\hline 
\end{tabular}
\par\end{centering}
\caption{Modal scales associated to a fixed base-chord}
\label{tab:modclass}
\end{table}

For instance, choose $[B]$ as a major seven chord, then it is possible to
associate to $[B]$ three different classes of tension-triads $[T_i]$.
Considering rows associated to the $maj7$ chord type in table
\ref{tab:modclass}, we have that the modes listed in the table are representable
as the following points in $\mathscr M$:
\begin{enumerate}
\item Ionian: $i:=([Cmaj7],[D-])$;
\item Lydian: $l:=([Cmaj7],[D])$;
\item Lydian $\sharp 2$: $l_{\sharp 2}:=([Cmaj7],[D\sharp -^{\flat 5}])$
\end{enumerate}
Every couple $([B],[T])$ identifies a point in $\mathscr M$ which can be
associated uniquely to a modal scale. The scale is given by the set of notes
$\left\{b_1,b_2,b_4,b_4,t_1,t_2,t_3\right\}$ where
$B=\left\{b_1,\dots,b_4\right\}$ and $T=\left\{t_1,t_2,t_3\right\}$. Following
\cite[Chapter 10]{Pis59} for the analysis of nonharmonic tones in classical
harmony we are entitled to define
\begin{defn}
Let $B=\left\{b_1,\dots,b_4\right\}$ and $T=\left\{t_1,t_2,t_3\right\}$. 
Non chord tones are notes which do not belong to the base-chord; i.e.
\[
t_i\in T \mbox{ such that } t_i\not\in B.
\]
\end{defn}

From definition \ref{def:modo} we know that for a certain mode $m\in\mathscr M$,
$m=([B],[T])$ it has to be $[B]\cap[T]=\emptyset$. So every note belonging to a
tension-triad associated to a base-chord is a non chord tone.
Every point in $\mathscr M$ has different musical meanings, here follows a list
of information we can read for each point $m\in\mathscr M$
\begin{enumerate}
\item $m$ identifies a unique mode;
\item chord tones and non chord tones are splitted into two components,
respectively $[B]$ and $[T]$
\item considering the notes belonging to $[B]$ and $[T]$ we deduce the modal
scale associated to $m$, that can be re-ordered in a $7$-uple in which the
degrees of the scale are displayed from the root, to the seventh note.
\end{enumerate}

Thus fixed a base-chord we can write down an ordered modal scale associated to
the chord for every available choice of tension-triad.

\begin{example}
Fix a seventh chord, for instance a Cmaj7. The idea is to split stable
and unstable notes of an heptatonic scale belonging to the \emph{base-chord},
as follows

\begin{table}[H]
\begin{centering}
\begin{tabular}{ccccccccccccc}
{\large C} & {\large $\rightarrow$} & {\large $\square$} & {\large
$\rightarrow$} & {\large E} & {\large $\rightarrow$} & {\large $\square$} 
& {\large $\rightarrow$} & {\large G} & {\large $\rightarrow$} & {\large
$\square$} & {\large $\rightarrow$} & {\large B}\tabularnewline
\end{tabular}
\par\end{centering}

\end{table}

\end{example}
White squares are placeholders for the note of a suitable triad. As we showed in the
previous paragraph, choosing a $D$ minor triad one can find the \emph{C ionian}
scale, considering a $D$ major triad we have a \emph{C lydian} and with a $D\sharp$ diminished triad we obtain the \emph{C lydian $\sharp 2$} scale.

\begin{table}[H]
\begin{centering}
{\large }%
\begin{tabular}{ccccccccccccc}
{\large C} & {\large $\rightarrow$} & {\large $ $D} & {\large $\rightarrow$} &
{\large E} & {\large $\rightarrow$}
& {\large F} & {\large $\rightarrow$} & {\large G} & {\large $\rightarrow$} &
{\large A} & {\large $\rightarrow$} & {\large B}\tabularnewline
\end{tabular}

\vspace{5mm}
\centering{}{\large }%
\begin{tabular}{ccccccccccccc}
{\large C} & {\large $\rightarrow$} & {\large D} & {\large $\rightarrow$} &
{\large E} & {\large $\rightarrow$} & {\large F$\sharp$} 
& {\large $\rightarrow$} & {\large G} & {\large $\rightarrow$} & {\large A} &
{\large $\rightarrow$} & {\large B}\tabularnewline
\end{tabular}
\par\end{centering}{\large \par}
\vspace{5mm}
\centering{}{\large }%
\begin{tabular}{ccccccccccccc}
{\large C} & {\large $\rightarrow$} & {\large D$\sharp$} & {\large
$\rightarrow$} & {\large E} & {\large $\rightarrow$} & {\large F$\sharp$} 
& {\large $\rightarrow$} & {\large G} & {\large $\rightarrow$} & {\large A} &
{\large $\rightarrow$} & {\large B}\tabularnewline
\end{tabular}
\end{table}


\section{A geometrical representation of modes through graphs}\label{sec:graphs}

This section ideally splits into two parts. In the first part, 
for the sake of the reader as well as for fixing our notations, we shall recall some 
well-known facts about graphs, maximal trees and fundamental group. For 
further details we shall refer to \cite{Gib10} (or any elementary textbook in algebraic topology) 
and references therein.

In the second part which is the core of the first part of the paper, we start 
defining the associated planar oriented graph to each base-chord. Since every mode 
is an heptatonic scale, the orientation of the graph is induced by the natural 
sequence of the degrees of the scale. 

\subsection{Some mathematical preliminaries}

\begin{defn}\label{def:grafoastratto}
An {\em abstract unoriented graph\/} is a pair $(V,E)$ where $V$ is a finite set
and $E$ is a set of 
unordered pairs of different elements of $V$. Thus an element of $E$ is of the
form $\{v,w\}$ where 
$v$ and $w$ belong to $V$ and $v \not=w$. We call vertices the 
elements of $V$ and edges the elements $\{v,w\}$ of $E$ connecting 
$v$ and $w$ (or $w$ and $v$).
\end{defn}
\begin{defn}\label{def:realizzazione}
Let $(V,E)$ be an abstract graph. A {\em realization\/} of $(V,E)$ is a set of
points in $\R^N$, one point for 
each vertex and segments joining precisely those pairs of
points which correspond to edges. The 
points are the {\em vertices\/} and the segments are the {\em edges\/}; the
realization is termed a {\em graph\/}. 
We require that the following two {\em intersection conditions\/} hold: 
\begin{enumerate}
 \item two edges meet either in a common end-point or at all;
 \item no vertex lies on an edge except at one of its ends. 
\end{enumerate}
We denote by $(vw)$ the edge joining $v$ and $w$.
\end{defn}
\begin{defn}\label{def:isomorfismo}
Two abstract (unoriented) graphs $(V,E)$ and $(V',E')$ are {\em
isomorphic\/} if there exists a bijective map $f:V \to V'$ 
such that 
\[
 \{v,w\} \in E \Longleftrightarrow \{f(v), f(w)\} \in E'.
\]
\end{defn}

\begin{rem}
Analogous definitions for oriented graphs are obtained by replacing
unordered pairs $\{\cdot, \cdot\}$ by 
ordered pairs $(\cdot, \cdot)$.
\end{rem}
\begin{defn}\label{def:pathinagraph}
For $n \geq 1$, a path on a graph $G$ from $v^1$ to $v^{n+1}$ is a sequence of
vertices and edges
 \[
  v^1e^1v^2e^2\dots v^n e^nv^{n+1}
 \]
where $e^1=(v^1v^2),$ $e^2=(v^2v^3),\dots, e^n=(v^n v^{n+1})$. 
\end{defn}
If $G$ is oriented, we only require that $e^i=(v^iv^{i+1})$ or
$e^i=(v^{i+1}v^i)$ for $i =1, \dots, n$; that is, the 
edges along the path are oriented in the opposite way.
\begin{defn}
The path is {\em simple\/} if  $e^1, \dots, e^n$ are all distinct, and
 $v^1, \dots, v^{n+1}$ are all 
distincts except that possibly $v^1= v^{n+1}$. If the simple path has
$v^1=v^{n+1}$ and $n>0$ is called a {\em loop\/}.

A graph $G$ is said to be {\em connected\/} if, given any two vertices $v$ and
$w$ of $G$ there is a path on $G$ from $v$ 
to $w$. A graph which is connected and without loops is called a {\em tree\/}.
\end{defn}
\begin{figure}[H]
\begin{centering}
\subfloat[An example of a planar grap]
{\begin{centering}
\includegraphics[width=3cm]{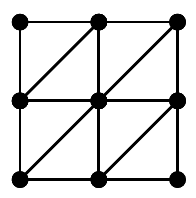}
\par\end{centering}
}$\;\;\;\;$\subfloat[A maximal tree]{\begin{centering}
\includegraphics[width=3cm]{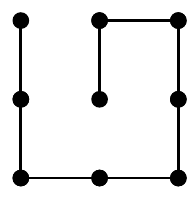}
\par\end{centering}
}
\par\end{centering}
\end{figure}
\begin{defn}
Given a graph  $G$, a graph $H$ is called a {\em subgraph\/} of $G$ is the
vertices of $H$ are vertices of $G$ and the 
edges of $H$ are edges of $G$. Also $H$ is called a {\em proper subgraph\/} of
$G$ if $H \not=G$. 
\end{defn}
The following definition is central in the sequel of the paper. 
\begin{defn}\label{def:admissible}
Let $G$ be a graph, $H$  be any maximal tree in $G$, $S$ be a subset of the 
vertices set and $k$ is an integer. An {\em admissible path $\gamma$ in $G$ with
respect to $S$ 
of  length $k$\/} is any  proper subgraph of $H$ satisfying the  two
conditions: 
\begin{enumerate}
 \item each vertex $v \in S$ lies in $\gamma$;
 \item the total number of vertices in $\gamma$ is $k$.
\end{enumerate}
\end{defn}
We also observe that any graph $G$ as a subgraph which is a tree (e.g. the empty
subgraph is a  tree) so that the set 
$\mathscr T$ of subgraphs of $G$ which are trees will have maximal elements.
That is, there exists at least one
$T \in \mathscr T$ such that $T$ is not a proper subgraph of any $T' \in
\mathscr T$.
\begin{lem}
Let $G$ be a connected graph. A subgraph $T$ of $G$ is a maximal tree for $G$ if
and only 
if $T$ is a tree containing 
all the vertices of $G$.
\end{lem}
\proof Cfr. \cite[Proposition 1.11, pag.18]{Gib10}.\finedim
For a connected graph $G$ there is a standard way to compute the homotopy group.
In fact the following result holds:
\begin{prop}\label{thm:sullomotopia}
For a connected graph $G$ with maximal tree $T$, $\pi_1(G)$ is a free group with
basis the classes $[f_\alpha]$ 
corresponding to the edges $e_\alpha$ of $X\setminus T$.
\end{prop}
\proof Cfr. \cite[Proposition 1A.2 pag.84]{Hat02}.\finedim 
Given a finite connected graph $G$, we denote by  $\chi(G)$ the {\em Euler characteristic\/} 
defined as the number of vertices minus the number of edges. From proposition 
\ref{thm:sullomotopia} we immediately get: 
\begin{cor}
Let $G$ be a finite connected graph. Then $\pi_1(G)$ is a free group over $1-\chi(G)$ generators.
More precisely the Euler characteristic of a finite connected graph is a topological invariant.
\end{cor}


\subsection{Graphs and base-chords}

\begin{defn}\label{def:associated_graph}
Given a base-chord $[B] \in \mathfrak C_4$ the associated graph $\mathscr G([B])$ is the
realization of the  abstract graph whose vertex set is given by  the set of
all notes forming $[B]$ and of every compatible 
tension-triad and the oriented edge set is represented by all possible oriented connections
between each vertex according to the order of the degrees of the scale; i.e. from the root to the seventh. (See figure \ref{fig:main}).

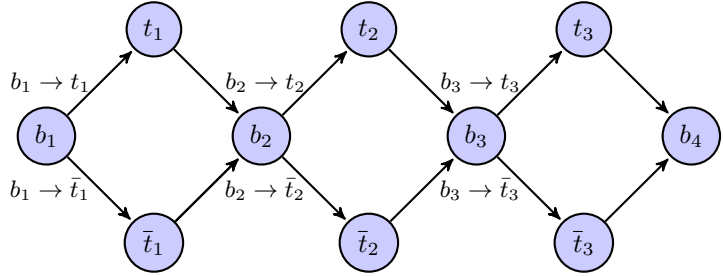
\begin{figure}[H]
\centering
\begin{tikzpicture}[->,>=stealth',shorten >=1pt,auto,node distance=2cm,
  thick,main node/.style={circle,fill=blue!20,draw,font=\sffamily\bfseries}]
  \node[main node] (1) {$b_1$};
  \node[main node] (2) [above right of=1] {$t_1$};
  \node[main node] (3) [below right of=1] {$\overline{t}_1$};
  \node[main node] (4) [below right of=2] {$b_2$};
  \node[main node] (5) [below right of=4] {$\overline{t}_2$};
 \node[main node] (6) [above right of=4] {$t_2$};
  \node[main node] (7) [below right of=6] {$b_3$};
  \node[main node] (8) [below right of=7] {$\overline{t}_3$};
  \node[main node] (9) [above right of=7] {$t_3$};
  \node[main node] (10) [above right of=8] {$b_4$};
  \path[every node/.style={font=\sffamily\small}]
    (1) edge node [left] {$b_1\rightarrow t_1$} (2)
    (1) edge node [left] {$b_1\rightarrow \bar{t}_1$} (3)
    (2) edge node [left] {} (4)
    (3) edge node [left] {} (4)
    (3) edge node [left] {} (4)
    (4) edge node [left] {$b_2\rightarrow \bar{t}_2$} (5)
    (4) edge node [left] {$b_2\rightarrow t_2$} (6)
    (6) edge node [left] {} (7)
    (5) edge node [left] {} (7)
    (7) edge node [left] {$b_3\rightarrow \bar{t}_3$} (8)
    (7) edge node [left] {$b_3\rightarrow t_3$} (9)
    (8) edge node [left] {} (10)
    (9) edge node [left] {} (10);        
\end{tikzpicture}
\caption{A graph built assuming (without loss of generality) that the modal
choices on a base-chord $B=\{b_1,b_2,b_3,b_4\}$ are given by two tension-triads
$T=\{t_1,t_2,t_3\}$ and
$\overline{T}=\{\overline{t}_1,\overline{t}_2,\overline{t}_3\}$}
\label{fig:main}
\end{figure}

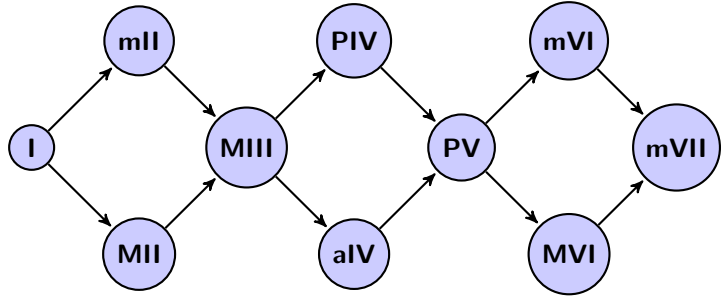
\begin{figure}[H]
\begin{tikzpicture}[->,>=stealth',shorten >=1pt,auto,node distance=2cm,
  thick,main node/.style={circle,fill=blue!20,draw,font=\sffamily\bfseries}]
  \node[main node] (1) {I};
  \node[main node] (2) [above right of=1] {mII};
  \node[main node] (3) [below right of=1] {MII};
  \node[main node] (4) [above right of=3] {MIII};
  \node[main node] (5) [above right of=4] {PIV};
  \node[main node] (6) [below right of=4] {aIV};
  \node[main node] (7) [above right of=6] {PV};
  \node[main node] (8) [above right of=7] {mVI};
  \node[main node] (9) [below right of=7] {MVI};
  \node[main node] (10) [above right of=9] {mVII};
  
  \path[every node/.style={font=\sffamily\small}]
    (1) edge node [left] {} (2)
    (1) edge node [left] {} (3)      
    (2) edge node [left] {} (4)      
    (3) edge node [left] {} (4)      
    (4) edge node [left] {} (5)      
    (4) edge node [left] {} (6)
    (5) edge node [left] {} (7)
    (6) edge node [left] {} (7)
    (7) edge node [left] {} (8)
    (7) edge node [left] {} (9)
    (8) edge node [left] {} (10)    
    (9) edge node [left] {} (10)    
            ;      
\end{tikzpicture}
\centering{}
\caption{A dominant chord graph (M:=major, m:=minor,  P:=perfect and a:=augmented).}
\end{figure}
\end{defn}
\begin{rem}
A priori the degrees can be arranged in the plane in infinitely many ways, giving rise to 
completely different unoriented graphs. However considering the orientation induced by the degrees 
of the scale all the oriented graphs are homeomorphic. Since homeomorphisms induce isomorphisms in 
homotopy, all of the homotopic classification as well as the topological measure of complexity, is 
not affected by the convention given in definition \ref{def:associated_graph}.

We also observe that on the second, fourth and sixth degrees we have at maximum two choices. This is a 
straightforward consequence of the constructions of the modes from the major, melodic minor and harmonic 
minor scales. (Read table \ref{tab:modclass} for further details).
\end{rem}

\begin{example} 
Labelling the previous graph using C as root note for each \emph{base-chord},
we obtain the following oriented planar graph.
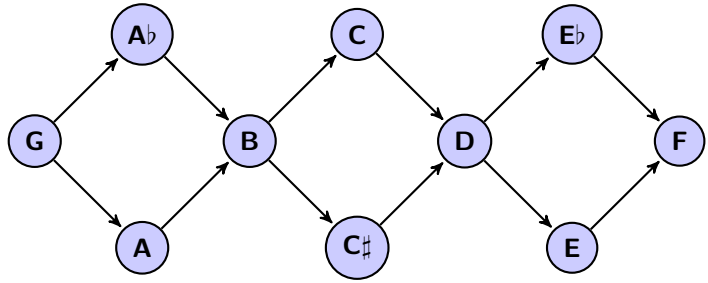
\begin{figure}[H]
\begin{tikzpicture}[->,>=stealth',shorten >=1pt,auto,node distance=2cm,
  thick,main node/.style={circle,fill=blue!20,draw,font=\sffamily\bfseries}]
  \node[main node] (1) {G};
  \node[main node] (2) [above right of=1] {A$\flat$};
  \node[main node] (3) [below right of=1] {A};
  \node[main node] (4) [above right of=3] {B};
  \node[main node] (5) [above right of=4] {C};
  \node[main node] (6) [below right of=4] {C$\sharp$};
  \node[main node] (7) [above right of=6] {D};
  \node[main node] (8) [above right of=7] {E$\flat$};
  \node[main node] (9) [below right of=7] {E};
  \node[main node] (10) [above right of=9] {F};
  
  \path[every node/.style={font=\sffamily\small}]
    (1) edge node [left] {} (2)
    (1) edge node [left] {} (3)      
    (2) edge node [left] {} (4)      
    (3) edge node [left] {} (4)      
    (4) edge node [left] {} (5)      
    (4) edge node [left] {} (6)
    (5) edge node [left] {} (7)
    (6) edge node [left] {} (7)
    (7) edge node [left] {} (8)
    (7) edge node [left] {} (9)
    (8) edge node [left] {} (10)    
    (9) edge node [left] {} (10)    
            ;      
\end{tikzpicture}
\centering{}
\caption{The graph associated to $G7$,
$\Gamma_{G7}$.}
\end{figure}

\end{example}

Following definition \ref{def:associated_graph} we can build a graph for each
type of chord\footnote{These types of chord are the one which are deduced from
the harmonization of the major, harmonic minor and melodic minor scale.}
\[
^{\circ 7},\,maj7^{\sharp5},\,-maj7,\,maj7,\,7,\,-7,\,-7^{\flat 5}.
\]
\begin{enumerate}
\item \textbf{Diminished seven: $\Gamma_{\circ7}$}. This kind of chord
appears only in the harmonization of the seventh degree of the harmonic
minor scale, therefore we have no choice except for the ultralocrian
mode. So, the graph associated to this seventh chord is 
\begin{figure}[H]
\begin{tikzpicture}[->,>=stealth',shorten >=1pt,auto,node distance=2cm,
  thick,main node/.style={circle,fill=blue!20,draw,font=\sffamily\bfseries}]
  \node[main node] (1) {I};
  \node[main node] (2) [right of=1] {mII};
  \node[main node] (3) [right of=2] {mIII};
  \node[main node] (4) [right of=3] {dIV};
  \node[main node] (5) [right of=4] {dV};
  \node[main node] (6) [right of=5] {mVI};
  \node[main node] (7) [right of=6] {dVII};
  
  \path[every node/.style={font=\sffamily\small}]
    (1) edge node [left] {} (2)
    (2) edge node [left] {} (3)      
    (3) edge node [left] {} (4)      
    (4) edge node [left] {} (5)      
    (5) edge node [left] {} (6)      
    (6) edge node [left] {} (7);      
\end{tikzpicture}
\centering{}\caption{The graph associated to diminished seventh chords,
$\Gamma_{\circ 7}$.}\label{fig: dim}
\end{figure}
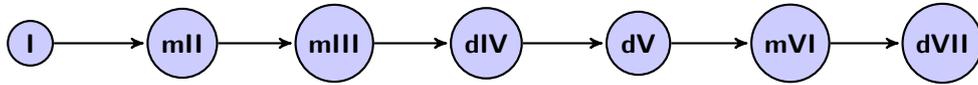
\begin{figure}[H]
\begin{tikzpicture}[->,>=stealth',shorten >=1pt,auto,node distance=2cm,
  thick,main node/.style={circle,fill=red!20,draw,font=\sffamily\bfseries}]
  \node[main node] (1) {I};
  \node[main node] (2) [right of=1] {mII};
  \node[main node] (3) [right of=2] {mIII};
  \node[main node] (4) [right of=3] {dIV};
  \node[main node] (5) [right of=4] {dV};
  \node[main node] (6) [right of=5] {mVI};
  \node[main node] (7) [right of=6] {dVII};

  \path[every node/.style={font=\sffamily\small}]
    (1) edge node [left] {} (2)
    (2) edge node [left] {} (3)      
    (3) edge node [left] {} (4)      
    (4) edge node [left] {} (5)      
    (5) edge node [left] {} (6)      
    (6) edge node [left] {} (7);      
\end{tikzpicture}
\centering{}\caption{A (unique) maximal tree for $\Gamma_{\circ 7}$.}
\end{figure}
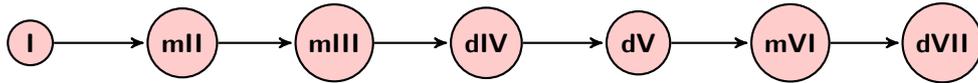

\item \textbf{Major seven $\sharp5$: $\Gamma_{{maj7^{\sharp5}}}$.} Fixing a
$maj7^{\,\sharp5}$ chord as \emph{base }of the mode,
we have two different possibilities: either the ionian sharp five
or the lydian sharp five modal scale. In terms of graph, we have 
\begin{figure}[H]

\begin{tikzpicture}[->,>=stealth',shorten >=1pt,auto,node distance=2cm,
  thick,main node/.style={circle,fill=blue!20,draw,font=\sffamily\bfseries}]

  \node[main node] (1) {I};
  \node[main node] (2) [right of=1] {MII};
  \node[main node] (3) [right of=2] {MIII};
  \node[main node] (4) [above right of=3] {aIV};
  \node[main node] (5) [below right of=3] {PIV};
  \node[main node] (6) [below right of=4] {aV};
  \node[main node] (7) [right of=6] {MVI};
  \node[main node] (8) [right of=7] {MVII};
  
  \path[every node/.style={font=\sffamily\small}]
    (1) edge node [left] {} (2)
    (2) edge node [left] {} (3)      
    (3) edge node [left] {} (4)
    (3) edge node [left] {} (5)      
    (4) edge node [left] {} (6)
    (5) edge node [left] {} (6)      
    (6) edge node [left] {} (7)      
    (7) edge node [left] {} (8);      
\end{tikzpicture}
\centering{}\caption{The graph associated to diminished seventh chords,
$\Gamma_{maj7^{\sharp 5}}$.}\label{fig: majd5}
\end{figure}
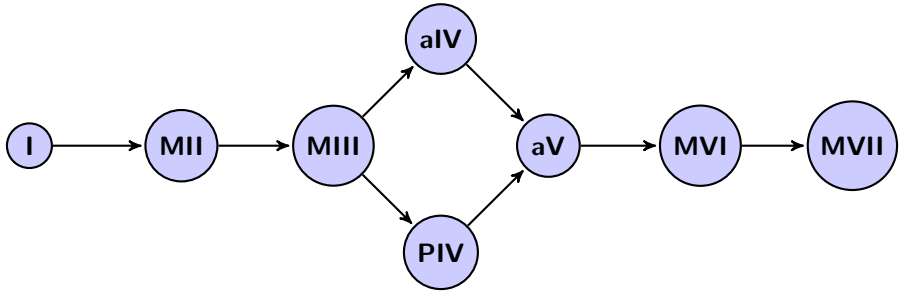
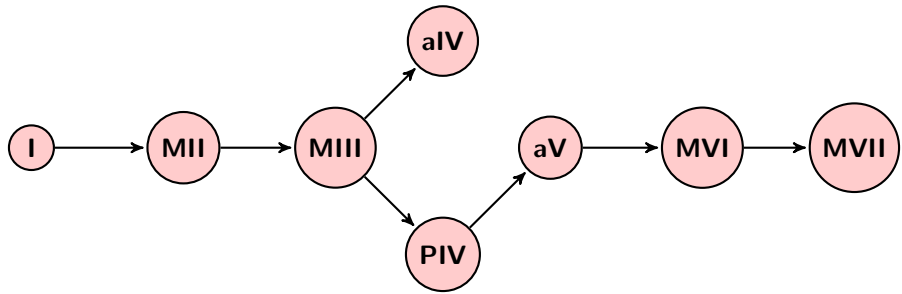
\begin{figure}[H]
\begin{tikzpicture}[->,>=stealth',shorten >=1pt,auto,node distance=2cm,
  thick,main node/.style={circle,fill=red!20,draw,font=\sffamily\bfseries}]
  \node[main node] (1) {I};
  \node[main node] (2) [right of=1] {MII};
  \node[main node] (3) [right of=2] {MIII};
  \node[main node] (4) [above right of=3] {aIV};
  \node[main node] (5) [below right of=3] {PIV};
  \node[main node] (6) [below right of=4] {aV};
  \node[main node] (7) [right of=6] {MVI};
  \node[main node] (8) [right of=7] {MVII};
  
  \path[every node/.style={font=\sffamily\small}]
    (1) edge node [left] {} (2)
    (2) edge node [left] {} (3)      
    (3) edge node [left] {} (4)
    (3) edge node [left] {} (5)      
    (5) edge node [left] {} (6)      
    (6) edge node [left] {} (7)      
    (7) edge node [left] {} (8);      
\end{tikzpicture}
\centering{}\caption{A maximal tree of $\Gamma_{maj7^{\sharp 5}}$.}
\end{figure}

\item \textbf{Minor major seven: $\Gamma_{-maj7}$}. In this case we
can choose between two different modes: hypoionian and hypoionian
$\flat6$. The graph is 
\begin{figure}[H]
\begin{tikzpicture}[->,>=stealth',shorten >=1pt,auto,node distance=2cm,
  thick,main node/.style={circle,fill=blue!20,draw,font=\sffamily\bfseries}]
  \node[main node] (1) {I};
  \node[main node] (2) [right of=1] {MII};
  \node[main node] (3) [right of=2] {mIII};
  \node[main node] (4) [right of=3] {PIV};
  \node[main node] (5) [right of=4] {PV};
  \node[main node] (6) [above right of=5] {MVI};
  \node[main node] (7) [below right of=5] {mVI};  
  \node[main node] (8) [below right of=6] {MVII};
  
  \path[every node/.style={font=\sffamily\small}]
    (1) edge node [left] {} (2)
    (2) edge node [left] {} (3)      
    (3) edge node [left] {} (4)      
    (4) edge node [left] {} (5)      
    (5) edge node [left] {} (6)
    (5) edge node [left] {} (7)
    (6) edge node [left] {} (8)      
    (7) edge node [left] {} (8);      
\end{tikzpicture}
\centering{}\caption{The graph associated to minor major seventh chords,
$\Gamma_{-maj 7}$.}\label{fig: mmaj7}
\end{figure}
\begin{figure}[H]
\begin{tikzpicture}[->,>=stealth',shorten >=1pt,auto,node distance=2cm,
  thick,main node/.style={circle,fill=red!20,draw,font=\sffamily\bfseries}]
  \node[main node] (1) {I};
  \node[main node] (2) [right of=1] {MII};
  \node[main node] (3) [right of=2] {mIII};
  \node[main node] (4) [right of=3] {PIV};
  \node[main node] (5) [right of=4] {PV};
  \node[main node] (6) [above right of=5] {MVI};
  \node[main node] (7) [below right of=5] {mVI};  
  \node[main node] (8) [below right of=6] {MVII};
  
  \path[every node/.style={font=\sffamily\small}]
    (1) edge node [left] {} (2)
    (2) edge node [left] {} (3)      
    (3) edge node [left] {} (4)      
    (4) edge node [left] {} (5)      
    (5) edge node [left] {} (6)
    (5) edge node [left] {} (7)
    (6) edge node [left] {} (8);      
\end{tikzpicture}
\centering{}\caption{A maximal tree associated to $\Gamma_{-maj 7}$.}
\end{figure}

\item \textbf{Major seven: $\Gamma_{maj7}$.} This is certainly a more
common chord than the previous ones. We expect to have more possibilities,
in fact a well known and simple \emph{base-chord} surely will bear
more \emph{tension-triads} than a naturally dissonant one. 
\begin{figure}[H]
\begin{tikzpicture}[->,>=stealth',shorten >=1pt,auto,node distance=2cm,
  thick,main node/.style={circle,fill=blue!20,draw,font=\sffamily\bfseries}]
  \node[main node] (1) {I};
  \node[main node] (2) [above right of=1] {aII};
  \node[main node] (3) [below right of=1] {MII};
  \node[main node] (4) [below right of=2] {MIII};
  \node[main node] (5) [above right of=4] {aIV};
  \node[main node] (6) [below right of=4] {PIV};
  \node[main node] (7) [above right of=6] {PV};
  \node[main node] (8) [right of=7] {MVI};
  \node[main node] (9) [right of=8] {MVII};
  
  \path[every node/.style={font=\sffamily\small}]
    (1) edge node [left] {} (2)
    (1) edge node [left] {} (3)      
    (2) edge node [left] {} (4)      
    (3) edge node [left] {} (4)          
    (4) edge node [left] {} (5)      
    (4) edge node [left] {} (6)
    (5) edge node [left] {} (7)      
    (6) edge node [left] {} (7)
    (7) edge node [left] {} (8)
    (8) edge node [left] {} (9);      
\end{tikzpicture}
\centering{}\caption{The graph associated to major seven chords,
$\Gamma_{maj7}$.}\label{fig: maj7}
\end{figure}
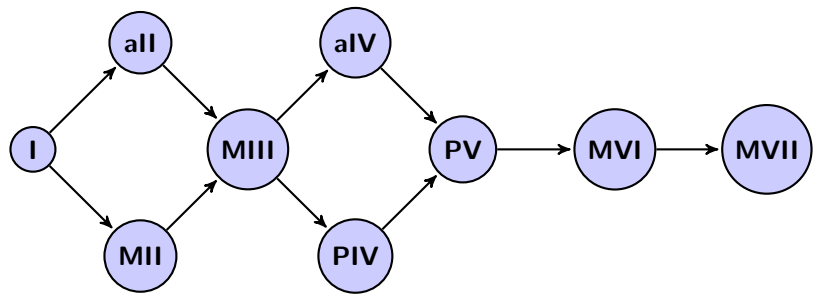
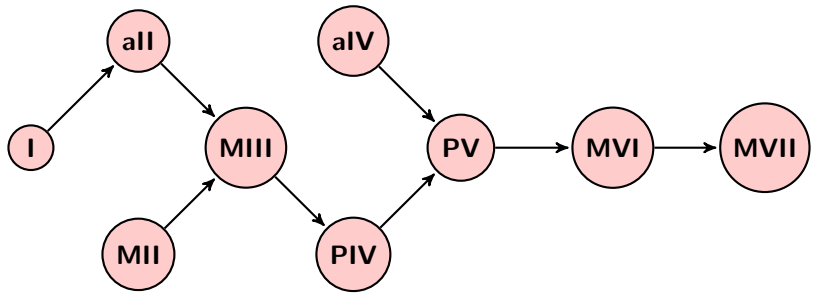
\begin{figure}[H]
\begin{tikzpicture}[->,>=stealth',shorten >=1pt,auto,node distance=2cm,
  thick,main node/.style={circle,fill=red!20,draw,font=\sffamily\bfseries}]
  \node[main node] (1) {I};
  \node[main node] (2) [above right of=1] {aII};
  \node[main node] (3) [below right of=1] {MII};
  \node[main node] (4) [below right of=2] {MIII};
  \node[main node] (5) [above right of=4] {aIV};
  \node[main node] (6) [below right of=4] {PIV};
  \node[main node] (7) [above right of=6] {PV};
  \node[main node] (8) [right of=7] {MVI};
  \node[main node] (9) [right of=8] {MVII};
  
  \path[every node/.style={font=\sffamily\small}]
    (1) edge node [left] {} (2)
    (2) edge node [left] {} (4)      
    (3) edge node [left] {} (4)          
    (4) edge node [left] {} (6)
    (5) edge node [left] {} (7)      
    (6) edge node [left] {} (7)
    (7) edge node [left] {} (8)
    (8) edge node [left] {} (9);      
\end{tikzpicture}
\centering{}\caption{A maximal tree associated to $\Gamma_{maj7}$.}
\end{figure}

\item \textbf{Dominant: $\Gamma_{7}$ .} Dominant chords are largely used
in blues and traditional jazz thanks to their capability of bearing
tensions. Thus we have:
\begin{figure}[H]
\begin{tikzpicture}[->,>=stealth',shorten >=1pt,auto,node distance=2cm,
  thick,main node/.style={circle,fill=blue!20,draw,font=\sffamily\bfseries}]
  \node[main node] (1) {I};
  \node[main node] (2) [above right of=1] {mII};
  \node[main node] (3) [below right of=1] {MII};
  \node[main node] (4) [above right of=3] {MIII};
  \node[main node] (5) [above right of=4] {aIV};
  \node[main node] (6) [below right of=4] {PIV};
  \node[main node] (7) [above right of=6] {PV};
  \node[main node] (8) [above right of=7] {mVI};
  \node[main node] (9) [below right of=7] {MVI};
  \node[main node] (10) [above right of=9] {mVII};
  
  \path[every node/.style={font=\sffamily\small}]
    (1) edge node [left] {} (2)
    (1) edge node [left] {} (3)      
    (2) edge node [left] {} (4)      
    (3) edge node [left] {} (4)      
    (4) edge node [left] {} (5)      
    (4) edge node [left] {} (6)
    (5) edge node [left] {} (7)
    (6) edge node [left] {} (7)
    (7) edge node [left] {} (8)
    (7) edge node [left] {} (9)
    (8) edge node [left] {} (10)    
    (9) edge node [left] {} (10)    
            ;      
\end{tikzpicture}
\centering{}\caption{The graph associated to dominant chords, $\Gamma_7$.}
\label{fig: 7}
\end{figure}
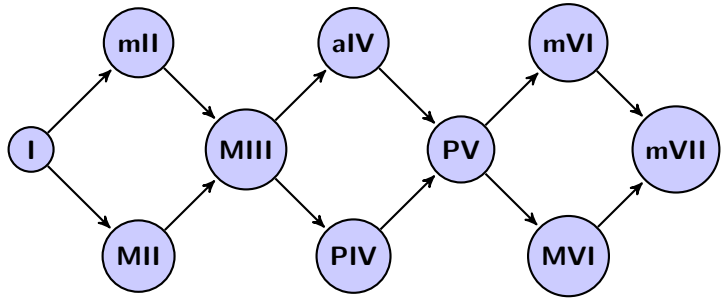
\begin{figure}[H]
\begin{tikzpicture}[->,>=stealth',shorten >=1pt,auto,node distance=2cm,
  thick,main node/.style={circle,fill=red!20,draw,font=\sffamily\bfseries}]
  \node[main node] (1) {I};
  \node[main node] (2) [above right of=1] {mII};
  \node[main node] (3) [below right of=1] {MII};
  \node[main node] (4) [above right of=3] {MIII};
  \node[main node] (5) [above right of=4] {aIV};
  \node[main node] (6) [below right of=4] {PIV};
  \node[main node] (7) [above right of=6] {PV};
  \node[main node] (8) [above right of=7] {mVI};
  \node[main node] (9) [below right of=7] {MVI};
  \node[main node] (10) [above right of=9] {mVII};
  
  \path[every node/.style={font=\sffamily\small}]
    (1) edge node [left] {} (2)
    (2) edge node [left] {} (4)      
    (3) edge node [left] {} (4)      
    (4) edge node [left] {} (6)
    (5) edge node [left] {} (7)
    (6) edge node [left] {} (7)
    (7) edge node [left] {} (8)
    (7) edge node [left] {} (9)
    (8) edge node [left] {} (10)    
            ;      
\end{tikzpicture}
\centering{}\caption{A maximal tree of $\Gamma_7$.}
\end{figure}
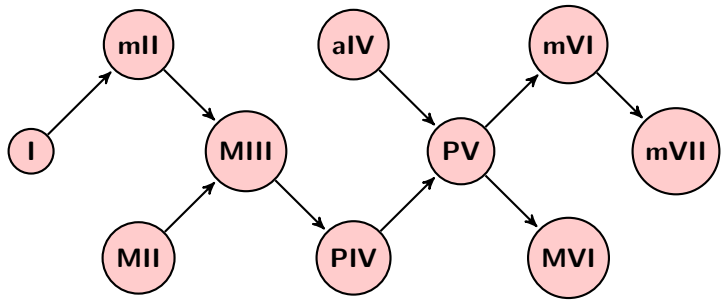

\item \textbf{Minor seven: $\Gamma_{-7}$.} For a minor seventh chord the only forbidden notes are 
the augmented second and the diminished fourth. So we have a graph isomorphic to $\Gamma_7$.
\begin{figure}[H]
\begin{tikzpicture}[->,>=stealth',shorten >=1pt,auto,node distance=2cm,
  thick,main node/.style={circle,fill=blue!20,draw,font=\sffamily\bfseries}]
  \node[main node] (1) {I};
  \node[main node] (2) [above right of=1] {mII};
  \node[main node] (3) [below right of=1] {MII};
  \node[main node] (4) [above right of=3] {mIII};
  \node[main node] (5) [above right of=4] {aIV};
  \node[main node] (6) [below right of=4] {PIV};
  \node[main node] (7) [above right of=6] {PV};
  \node[main node] (8) [above right of=7] {mVI};
  \node[main node] (9) [below right of=7] {MVI};
  \node[main node] (10) [above right of=9] {mVII};
  
  \path[every node/.style={font=\sffamily\small}]
    (1) edge node [left] {} (2)
    (1) edge node [left] {} (3)      
    (2) edge node [left] {} (4)      
    (3) edge node [left] {} (4)      
    (4) edge node [left] {} (5)      
    (4) edge node [left] {} (6)
    (5) edge node [left] {} (7)
    (6) edge node [left] {} (7)
    (7) edge node [left] {} (8)
    (7) edge node [left] {} (9)
    (8) edge node [left] {} (10)    
    (9) edge node [left] {} (10)    
            ;      
\end{tikzpicture}
\centering{}\caption{The graph associated to minor seven chords, $\Gamma_{-7}$.}
\label{fig: m7}
\end{figure}
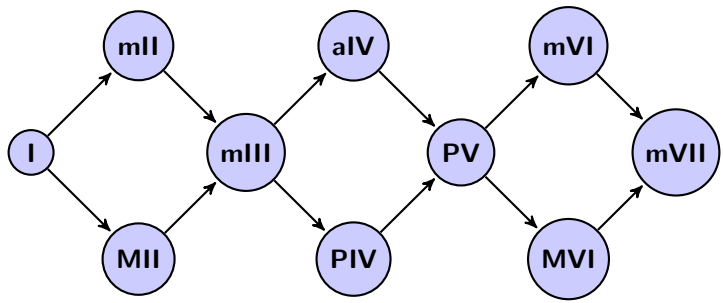
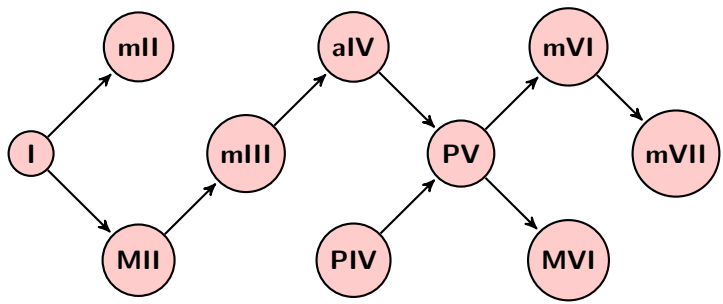
\begin{figure}[H]
\begin{tikzpicture}[->,>=stealth',shorten >=1pt,auto,node distance=2cm,
  thick,main node/.style={circle,fill=red!20,draw,font=\sffamily\bfseries}]
  \node[main node] (1) {I};
  \node[main node] (2) [above right of=1] {mII};
  \node[main node] (3) [below right of=1] {MII};
  \node[main node] (4) [above right of=3] {mIII};
  \node[main node] (5) [above right of=4] {aIV};
  \node[main node] (6) [below right of=4] {PIV};
  \node[main node] (7) [above right of=6] {PV};
  \node[main node] (8) [above right of=7] {mVI};
  \node[main node] (9) [below right of=7] {MVI};
  \node[main node] (10) [above right of=9] {mVII};
  
  \path[every node/.style={font=\sffamily\small}]
    (1) edge node [left] {} (2)
    (1) edge node [left] {} (3)      
    (3) edge node [left] {} (4)      
    (4) edge node [left] {} (5)      
    (5) edge node [left] {} (7)
    (6) edge node [left] {} (7)
    (7) edge node [left] {} (8)
    (7) edge node [left] {} (9)
    (8) edge node [left] {} (10)    
            ;      
\end{tikzpicture}
\centering{}\caption{A maximal tree of $\Gamma_{-7}$.}
\end{figure}
\item \textbf{Minor seven $\flat5$: $\Gamma_{-7\flat5}$.} 
In this case the \emph{base-chord} makes the difference: the root note and the diminished fifth form a tritone interval which gives a \emph{stable} sense of dissonance to the half-diminished seventh chord, that is emphasized by the minor second which
is natural in three of the four modal solutions we find on this type of chord, i.e. locrian, superlocrian and locrian $\sharp6$ scales. 
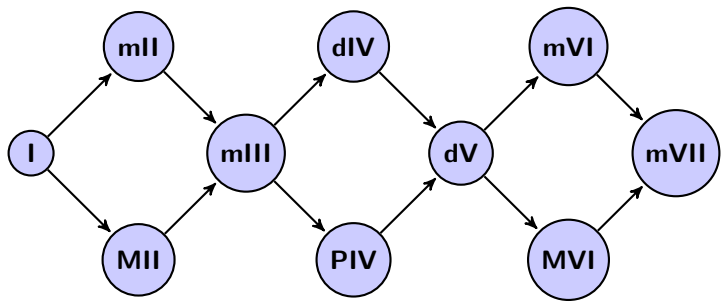
\begin{figure}[H]
\begin{tikzpicture}[->,>=stealth',shorten >=1pt,auto,node distance=2cm,
  thick,main node/.style={circle,fill=blue!20,draw,font=\sffamily\bfseries}]
  \node[main node] (1) {I};
  \node[main node] (2) [above right of=1] {mII};
  \node[main node] (3) [below right of=1] {MII};
  \node[main node] (4) [above right of=3] {mIII};
  \node[main node] (5) [above right of=4] {dIV};
  \node[main node] (6) [below right of=4] {PIV};
  \node[main node] (7) [above right of=6] {dV};
  \node[main node] (8) [above right of=7] {mVI};
  \node[main node] (9) [below right of=7] {MVI};
  \node[main node] (10) [above right of=9] {mVII};
  
  \path[every node/.style={font=\sffamily\small}]
    (1) edge node [left] {} (2)
    (1) edge node [left] {} (3)      
    (2) edge node [left] {} (4)      
    (3) edge node [left] {} (4)      
    (4) edge node [left] {} (5)      
    (4) edge node [left] {} (6)
    (5) edge node [left] {} (7)
    (6) edge node [left] {} (7)
    (7) edge node [left] {} (8)
    (7) edge node [left] {} (9)
    (8) edge node [left] {} (10)    
    (9) edge node [left] {} (10)    
            ;      
\end{tikzpicture}
\centering{}\caption{The graph associated to minor seven flat five chords,
$\Gamma_{-7^{\flat 5}}$.}\label{fig: m7b5}
\end{figure}
\begin{figure}[H]
\begin{tikzpicture}[->,>=stealth',shorten >=1pt,auto,node distance=2cm,
  thick,main node/.style={circle,fill=red!20,draw,font=\sffamily\bfseries}]
  \node[main node] (1) {I};
  \node[main node] (2) [above right of=1] {mII};
  \node[main node] (3) [below right of=1] {MII};
  \node[main node] (4) [above right of=3] {mIII};
  \node[main node] (5) [above right of=4] {dIV};
  \node[main node] (6) [below right of=4] {PIV};
  \node[main node] (7) [above right of=6] {dV};
  \node[main node] (8) [above right of=7] {mVI};
  \node[main node] (9) [below right of=7] {MVI};
  \node[main node] (10) [above right of=9] {mVII};
  
  \path[every node/.style={font=\sffamily\small}]
    (1) edge node [left] {} (2)
    (1) edge node [left] {} (3)      
    (3) edge node [left] {} (4)      
    (4) edge node [left] {} (5)      
    (5) edge node [left] {} (7)
    (6) edge node [left] {} (7)
    (7) edge node [left] {} (8)
    (8) edge node [left] {} (10)    
    (9) edge node [left] {} (10)    
            ;      
\end{tikzpicture}
\centering{}\caption{A maximal tree of $\Gamma_{-7^{\flat 5}}$.}
\end{figure}
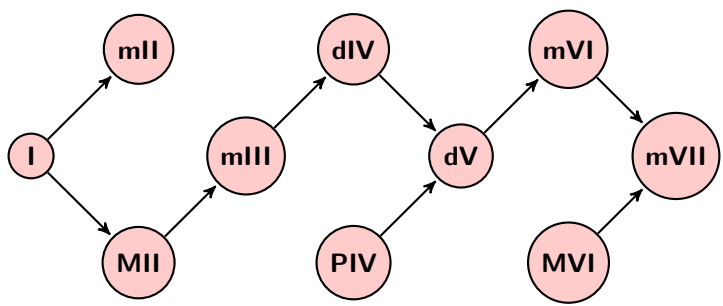
\end{enumerate}
\begin{example}
In this example we display three of the graphs investigated above, labelling the vertices with notes 
instead of the degrees of the modal scale.
\begin{figure}[H]
\begin{centering}
\subfloat[Minor major seven modes - Given a $-maj7$ chord, available modes
are ipoionic and ipoionic $\flat6$.]
{\begin{centering}
\includegraphics[width=5cm]{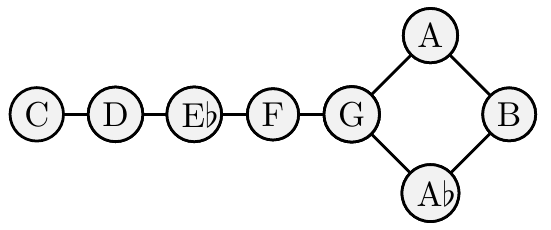}
\par\end{centering}
}$\;\;\;\;$\subfloat[Major seven modes - This graph represent modes on a major
seven
chords, $Cmaj7$ in this case.]{\begin{centering}
\includegraphics[width=5cm]{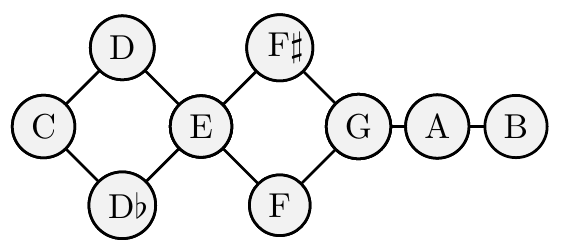}
\par\end{centering}
}
\par\end{centering}
\centering{}\subfloat[Ultra locrian - In this example tension-triad's notes are
double circled.
The base-chord is
$C^{\,\circ7}$.]{\centering{}\includegraphics[width=6cm]{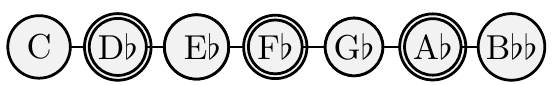}}
\end{figure}
\end{example}
\begin{rem}
It is clear by previous discussion that even if the graphs associated to 
$\Gamma_{7}$, $\Gamma_{-7}$, $\Gamma_{-7^{\,\flat5}}$ are isomorphic,  
they are built on different notes and hence they are quite different in the essence so
the homotopy is not suitable  to distinguish among them.
\end{rem}
All of these graphs show in a clear and net way how to construct new
modes from the 
existing ones. Given a graph $G$, let us consider any proper tree $H$ (not
maximal, in general!) 
contained  in $G$ and having $7$ vertices. By taking into account definition
\ref{def:admissible}, 
we give the following:
\begin{defn}\label{def:admissible-path}
Let $[B]$ a base-chord and $\mathscr G([B])$ be the associated base-chord graph.
An 
{\em admissible mode\/} is any admissible connected subgraph (or path in the graph) 
$\gamma([B])$ in $\mathscr G$ with respect to $[B]$ 
of length $7$. If $\gamma([B])$ is not a mode constructed above, we refer to as
{\em admissible special mode.\/}
\end{defn}
 Summing up, the previous discussion, we get:
\begin{prop}\label{prop:specmodes}
Given the base-chord $[B] \in \mathfrak C_4$ the following modes are the only 
admissible special modes:
\begin{enumerate}
\item if $[\Gamma_B]=[\Gamma_{maj7}]$ then $\gamma([B])$ is the path
 \[
 \gamma_{Ion\sharp 2}:=\{I,aII,MIII,PIV,PV,MVI,MVII\}
 \]
 \item if $[\Gamma_B]=[\Gamma_7]$ then $\gamma([B])$ are the paths
\begin{eqnarray*}
   \gamma_{Mix\flat 2} &:=&\{I,mII,MIII,PIV,PV,MVI,mVII\}\\
   \gamma_{Mix\flat 2\sharp 4}&:=&\{I,mII,MIII,aIV,PV,MVI,mVII\}\\
   \gamma_{Mix\sharp 4\flat 6}&:=&\{I,MII,MIII,aIV,PV,mVI,mVII\}\\
   \gamma_{Mix\flat 2\sharp 4\flat 6}&:=&\{I,mII,MIII,aIV,PV,mVI,mVII\}
\end{eqnarray*} 
 
 \item if $[\Gamma_B]=[\Gamma_{-7}]$ then $\gamma([B])$ are the paths
 \begin{eqnarray*}
   \gamma_{Eol\flat 2}&:=&\{I,mII,mIII,PIV,PV,mVI,mVII\}\\
   \gamma_{Eol\sharp 4}&:=&\{I,MII,mIII,aIV,PV,mVI,mVII\}\\
   \gamma_{Phr\sharp 4}&:=&\{I,mII,mIII,aIV,PV,mVI,mVII\}
 \end{eqnarray*}
 \item if $[\Gamma_B]=[\Gamma_{-7\flat 5}]$ then $\gamma([B])$ are the paths
 \begin{eqnarray*}
   \gamma_{Loc\sharp 2\sharp 6}&:=&\{I,mII,mIII,PIV,dV,mVI,mVII\}\\
   \gamma_{Sup\sharp 2}&:=&\{I,MII,mIII,dIV,dV,mVI,mVII\}\\
   \gamma_{Sup\sharp 6}&:=&\{I,mII,mIII,dIV,dV,MVI,mVII\}\\
   \gamma_{Sup\sharp 2\sharp 6}&:=&\{I,MII,mIII,dIV,dV,MVI,mVII\}
 \end{eqnarray*} 
\end{enumerate}
\end{prop}
\proof The result readily follows by the previous graph classification. Let us
consider every class of chord to prove the existence of special modes.
\begin{itemize}
\item $^{\circ 7}$. It is not possible to have path associated to special modes on
the graph $\Gamma_{\circ 7}$ (figure \ref{fig: dim}), since there is only one
path available which represent the ultra locrian mode.
\item $maj7^{\sharp 5}$ and $-maj7$. In both $\Gamma_{maj7^{\sharp 5}}$ (figure
\ref{fig: majd5}) and $\Gamma_{-maj7}$ (figure \ref{fig: mmaj7}) The only
available choice is on the fourth and the sixth degree of the modal scale,
respectively. This fact implies that only two admissible modes can be built on
such graph and they differs exactly for one note. So we can choose among two
paths on the graph which are exactly the two admissible modes we used to build
the graph.
\item $maj7$. In $\Gamma_{maj7}$ (figure \ref{fig: maj7}) there are $2^2$
available choices. The modes which generate this graph are three, so there is a
special mode which represent the admissible path on $\Gamma_{maj7}$ which is
different from the paths representing the Ionian, Lydian and Lydian $\sharp 2$
scales. The only possible, admissible path is 
\[
\{I,aII,MIII,PIV,PV,MVI,MVII\}.
\]
\item $7$. Four admissible non special modes generate $\Gamma_{7}$ (see table 
\ref{tab:modclass} and figure \ref{fig: 7}). The total number of admissible
modes in this graph is $2^3=8$. We expect to find $4$ special modes:
\begin{eqnarray*}
\{I,mII,MIII,PIV,PV,MVI,mVII\}\\
\{I,mII,MIII,aIV,PV,MVI,mVII\}\\
\{I,MII,MIII,aIV,PV,mVI,mVII\}\\
\{I,mII,MIII,aIV,PV,mVI,mVII\}
\end{eqnarray*} 
\item $-7$ and $-7^{\flat 5}$. This cases are similar to the previous one.
$\Gamma_{-7}$ (figure \ref{fig: m7}) is generated by $5$ admissible, non special
modes (table \ref{tab:modclass}), so we have $3$ special modes which are
\begin{eqnarray*}
\{I,mII,mIII,PIV,PV,mVI,mVII\}\\
\{I,MII,mIII,aIV,PV,mVI,mVII\}\\
\{I,mII,mIII,aIV,PV,mVI,mVII\}.
 \end{eqnarray*}
$\Gamma_{-7^{\flat 5}}$ (figure \ref{fig: m7b5}) is generated by four admissible non
special modes (table \ref{tab:modclass}), we have the following four special
modes:
 \begin{eqnarray*}
\{I,mII,mIII,PIV,dV,mVI,mVII\}\\
\{I,MII,mIII,dIV,dV,mVI,mVII\}\\
\{I,mII,mIII,dIV,dV,MVI,mVII\}\\
\{I,MII,mIII,dIV,dV,MVI,mVII\}
 \end{eqnarray*} 
\end{itemize}
\finedim
\subsection{Topological complexity measure of a base-chord}

The aim of this section is to associate a measure of complexity to each base
chord. 

\begin{defn}\label{def:tipotopologico}
Given a base-chord $[B]$, let $\mathscr G([B])$ be the associated base-chord graph. We
call   {\em topological 
complexity measure \/} (t.c.m., for brevity) of $[B]$, i.e. $\tau([B])$, 
the number of generators of the fundamental group  of $\mathscr G([B])$.
\end{defn}
\begin{lem}
Let $[B]$ be a base-chord. The integer $\tau([B])$ is well-defined.
\end{lem}
\proof It is enough to observe that  given any base-chord $[B]$,  
the associated integer $\tau([B])$ is uniquely defined. In fact 
by the classification given in section \ref{sec:graphs},  at each base-chord
$[B]$ we can  uniquely associate a planar connected graph $\mathscr G([B])$. As direct consequence 
of proposition \ref{thm:sullomotopia}, the fundamental group of
$\pi_1\big(\mathscr G([B])\big)$ is a 
free group having $\tau([B])$ generators. 
\finedim 
\begin{prop}\label{prop:calcoloomotoopia}
Let $[B]$ be a base-chord, $\mathscr G([B])$ be a planar graph. Then the fundamental
group 
$\pi_1\big(\mathscr G([B])\big)$ and the topological measure of complexity are given
below:
\begin{table}[H]
\centering{}%
\begin{tabular}{ccc}
\toprule 
Base-chord [B]& $\pi_{1}\big([B]\big)$ & $\tau([B])$\tabularnewline
\midrule
\midrule 
$^{\circ7}$ & $\left\{ 1\right\} $ & 0\tabularnewline
\midrule 
$maj7^{\sharp5}$ & $\mathbb{Z}$ & 1\tabularnewline
\midrule 
$-maj7$ & $\mathbb{Z}$ & 1\tabularnewline
\midrule 
$maj7$ & $\mathbb{Z}^{\ast 2}$ & 2\tabularnewline
\midrule 
$7$ & $\mathbb{Z}^{\ast 3}$ & 3\tabularnewline
\midrule 
$-7$ & $\mathbb{Z}^{\ast 3}$ & 3\tabularnewline
\midrule 
$-7^{\flat5}$ & $\mathbb{Z}^{\ast 3}$ & 3\tabularnewline
\bottomrule
\end{tabular}
\end{table}
\end{prop}
\proof The proof follows from the classification given in section
\ref{sec:graphs} and 
proposition \ref{thm:sullomotopia}.\finedim
By recalling a classical result from algebraic topology, namely that collapsing
a contractible 
subcomplex is a homotopy equivalence, it is possible to provide another proof of
the proposition 
\ref{thm:sullomotopia}. In fact, the following holds:
\begin{lem}\label{thm:lemmafava}
 If the pair $(X,A)$ satisfies the homotopy extension property and $A$ is
contractible, then the 
 quotient map $q: X \to X/A$ is a homotopy equivalence.
\end{lem}
For  the proof see, for instance \cite[Propoisition 0.17, pag.15]{Hat02}.
Applying lemma \ref{thm:lemmafava} to the pair $(\mathscr G([B]),T([B]))$ where $T([B])$
is a maximal 
tree contained in the graph $\mathscr G([B])$, the first homotopy group of the graph
$\mathscr G([B])$ can be 
easily computed (by using the Seifert Van-Kampen theorem) by observing that
the 
projection to the quotient: 
\[
 \Pi: \mathscr G([B]) \to \mathscr G([B])/T([B])
\]
is an homotopy equivalence and that the quotient is homotopic to the wedge of
$\tau([B])$ times 
$\mathbb S^1$  we get the result. 


\section{Harmonic progressions and braids}\label{sec:braids}

Here we start recalling some definitions from braid theory, then we give a braids theoretical interpretation 
of modes constructed over a fixed base-chord and by  
using the concatenation property for braids we
analyse some chord progressions. Our basic reference is \cite{Han89}.

\subsection{Braids and links}

\begin{defn}
 A {\em braid\/} $\beta$ on $n$ strands is a collection of embeddings 
 \[
  \mathscr B:=\{\beta^\alpha; \beta^\alpha:[0,1] \to \R^3, \ \alpha =1, \dots,
n\}
 \]
with disjoint images such that: 
\begin{itemize}
 \item $\beta^\alpha(0)=(0,\alpha,0)$;
 \item $\beta^\alpha(1)=(1, \tau(\alpha),0)$ for some permutation $\tau$;
 \item the images of each $\beta^\alpha$ is transverse to all planes
$\{x=const\}$, 
 where $\R^3$ is equipped with the $(x,y,z)$ coordinates.
\end{itemize}
\end{defn}
\begin{figure}[H]
\begin{centering}
\subfloat[$\sigma_1$]
{\begin{centering}
\includegraphics[height=3cm]{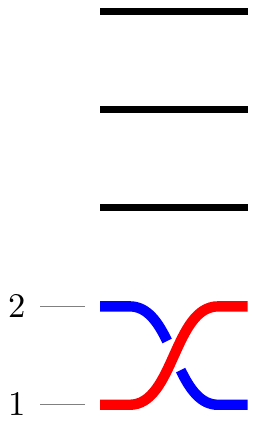}
\par\end{centering}
}$\;\;\;\;$\subfloat[$\sigma_2$]{\begin{centering}
\includegraphics[height=3cm]{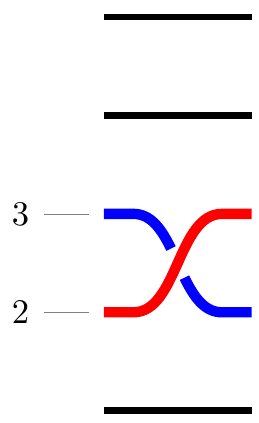}
\par\end{centering}
}$\;\;\;\;$\subfloat[$\sigma_i$]{\begin{centering}
\includegraphics[height=3cm]{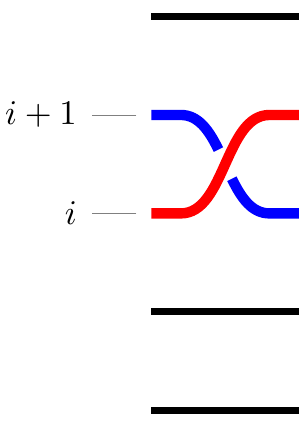}
\par\end{centering}
}
\par 
\end{centering}
\caption{The action of generators $\{\sigma_i\}_{i=1}^n$ on the $i^{th}$
strand.}
\label{fig:braidgen}
\end{figure}
\begin{defn}
Two such braids are said to lie in the same {\em topological braid class\/} if
they 
are homotopic in the sense of braids: one can deform one braid to the other
without any 
intersection among the strands. 
\end{defn}
There is a natural group structure on the space of topological braids with $n$
strands, $B_n$, 
given by concatenation. Using generators $\sigma_i$ which interchanges the
$i-th$ and $(i+1)-th$ strands 
with a positive crossing yields the presentation for $B_n$. 
\begin{figure}[H]
\begin{centering}
\subfloat[First braid]
{\begin{centering}
\includegraphics[height=2cm]{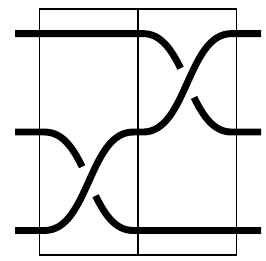}
\par\end{centering}
}$\;\;\;\;$\subfloat[Second braid]{\begin{centering}
\includegraphics[height=2cm]{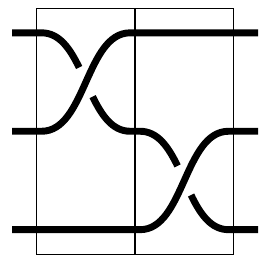}
\par\end{centering}
}
$\;\;\;\;$\subfloat[Their concatenation]{\begin{centering}
\includegraphics[height=2cm]{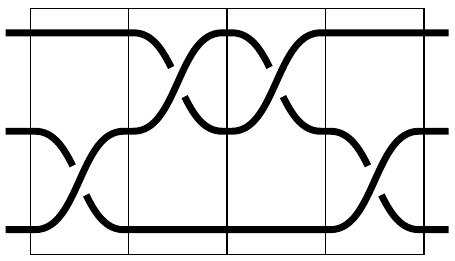}
\par\end{centering}
}
\caption{Concatenation of braids}
\par\end{centering}
\end{figure}
We let 
\[
\begin{split}
p_1:& \quad \sigma_i\,\sigma_j=\sigma_j\,\sigma_i; \quad |i-j|>1\\
p_2:& \quad \sigma_i\, \sigma_{i+1} \,\sigma_i = \sigma_{i+1}\, \sigma_i \,
\sigma_{i+1}; \quad i<n-1.
\end{split}
\]
From a graphical point of view, the first property can be visualized as
follows: 
\begin{figure}[H]
\begin{centering}
\subfloat[$p_1$ property: $\sigma_i \,\sigma_j$]
{\begin{centering}
\includegraphics[height=3cm]{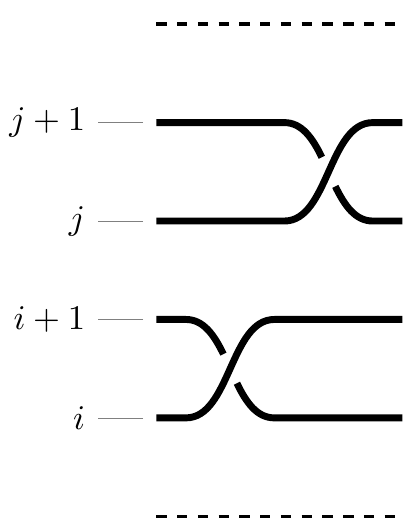}
\par\end{centering}
}$\;\;\;\;$\subfloat[$p_1$ property: $\sigma_j\, \sigma_i$ ]{\begin{centering}
\includegraphics[height=3cm]{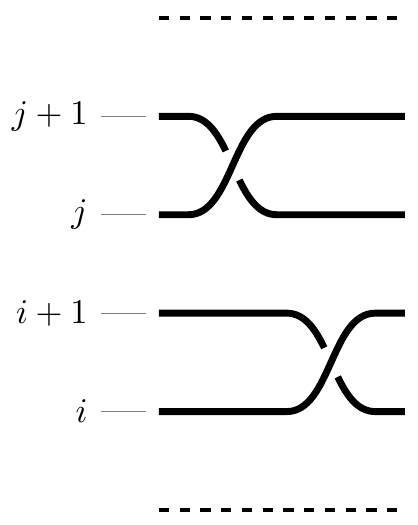}
\par\end{centering}
}
\par\end{centering}
\end{figure}
The second properties can be displayed as follows: 
\begin{figure}[H]
\begin{centering}
\subfloat[$p_2$ property: $\sigma_i\sigma_{i+1}\sigma_i$]
{\begin{centering}
\includegraphics[width=3cm]{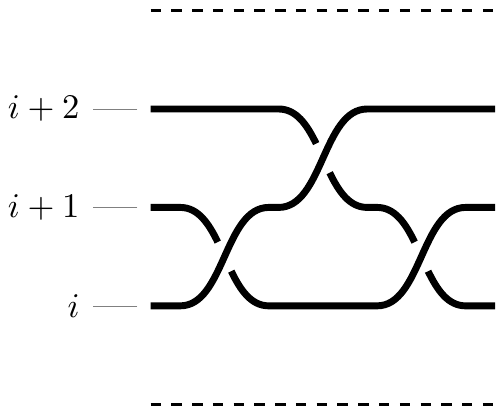}
\par\end{centering}
}$\;\;\;\;$\subfloat[$p_2$ property: $\sigma_{i+1}\sigma_{i}\sigma_{i+1}$]
{\begin{centering}
\includegraphics[width=3cm]{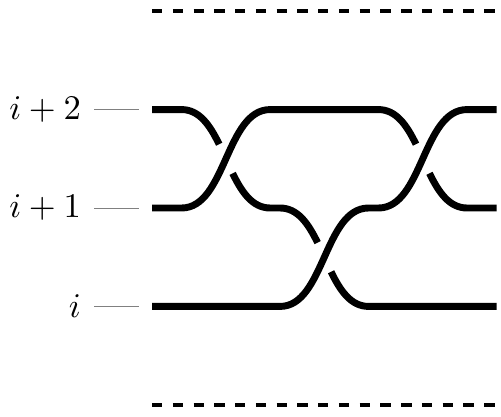}
\par\end{centering}
}
\par\end{centering}
\end{figure}
Then the presentation of $B_n$ is as follows:
\[
 B_n:=\langle  \sigma_1, \dots, \sigma_{n-1}: \  p_1 \ \textrm{and}\  p_2 \
\textrm{hold}\rangle.
\]
\begin{figure}[H]
\begin{centering}
\subfloat[The generator $\sigma_i$]
{\begin{centering}
\includegraphics[height=5cm]{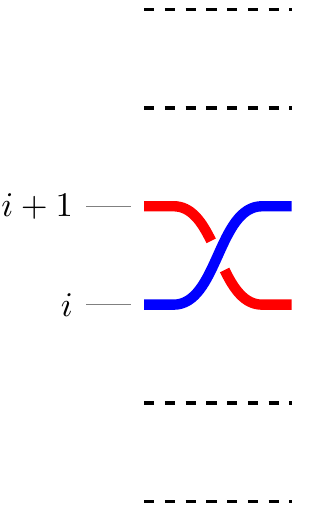}
\par\end{centering}
}$\;\;\;\;$\subfloat[The inverse $\sigma_i^{-1}$]{\begin{centering}
\includegraphics[height=5cm]{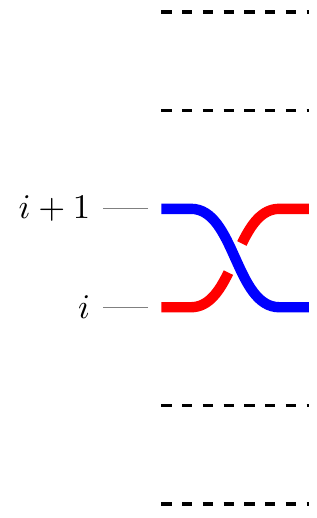}
\par\end{centering}
}
\par\end{centering}
\end{figure}

Now we consider $n$  distinct points  in the plane, i.e. $(x_j)_{j=1}^n\subset
\R^2$. Let
\[
 \Delta_{ij}=\{x=(x_1, \dots, x_n) \in \R^{2n}: x_i = x_j\}, \qquad i \not= j
\]
and borrowing the terminology of celestial mechanics, we can define the {\em
collision set\/} as: 
\[
\Delta= \bigcup_{i,j=1}^n \Delta_{ij}.
\]
Its complement $\widehat \chi_n(\R^2):= \R^{2n}\setminus \Delta$ is called the
{\em (collision-free) 
configuration space\/}. It can be thought as the configuration space for a set
of $n$ points without 
collisions. We equip $\widehat\chi_n(\R^2)$ with the topology induced from the
topology 
of the Euclidean space. Since the configuration space is the complement of a
finite 
union of two dimensional linear subspaces, by dimensional arguments, 
readily follows that $\widehat\chi_n(\R^2)$ is connected. 

There is a natural right action of $\mathfrak S(n)$ on the configuration space
$\widehat \chi_n(\R^2)$:
\[
 \mu: \widehat \chi_n(\R^2) \times \mathfrak S(n) \to \widehat \chi_n(\R^2) 
\]
defined by permutation of coordinates, i.e. 
\[
 \mu((x_1, x_n), \sigma)=(x_1, \dots, x_n)\cdot \sigma =(x_{\sigma(1)}, \dots,
x_{\sigma(n)}).
\]
It is to prove that $\mu$ is indeed a right and free action. Denote the orbit
space for the 
free action $\mu$ by $C_n(\R^2)$, or in standard notation, 
\[
 C_n(\R^2)= \widehat \chi_n(\R^2)/\mathfrak S(n).
\]

\begin{prop}
The (Artin) braid group $B_n$ can be canonically identified with the fundamental
group 
\[
 \pi_1(C_n(\R^2),c_0).
\]
\end{prop}
\proof Cfr. \cite[Theorem 3.1, pag 19]{Han89}.\finedim

We close the paragraph showing a braid theoretical representation of the planar oriented graphs 
associated to the base-chords of section \ref{sec:grafi}.

\begin{figure}[H]
\centering
\subfloat[$C^{\circ 7}$]
{\includegraphics[height=5cm]{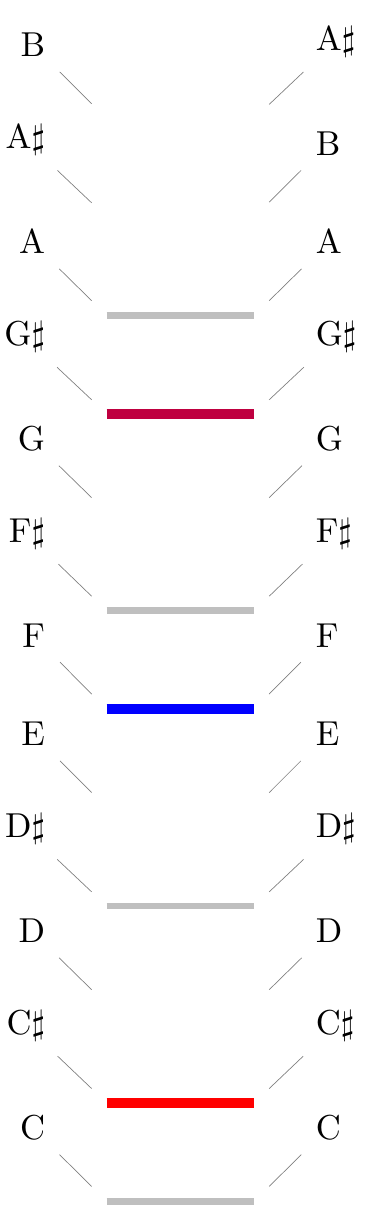}}
\hspace{15mm}
\subfloat[$Cmaj7^{\sharp 5}$]
{\includegraphics[height=5cm]{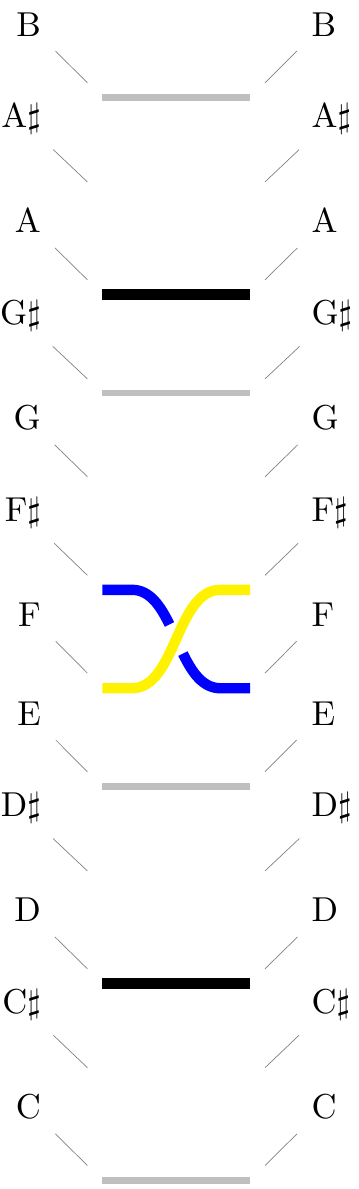}}
\hspace{15mm}
\subfloat[$C-maj7$]
{\includegraphics[height=5cm]{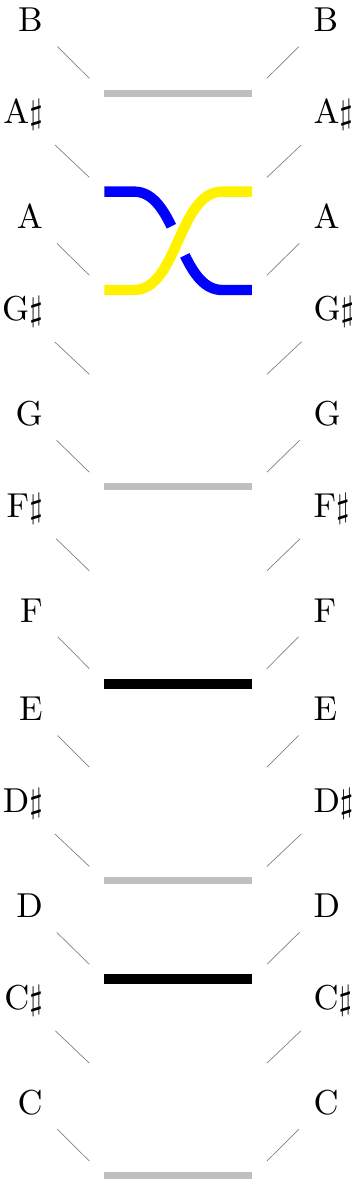}}
\hspace{15mm}
\subfloat[$Cmaj7$]
{\includegraphics[height=5cm]{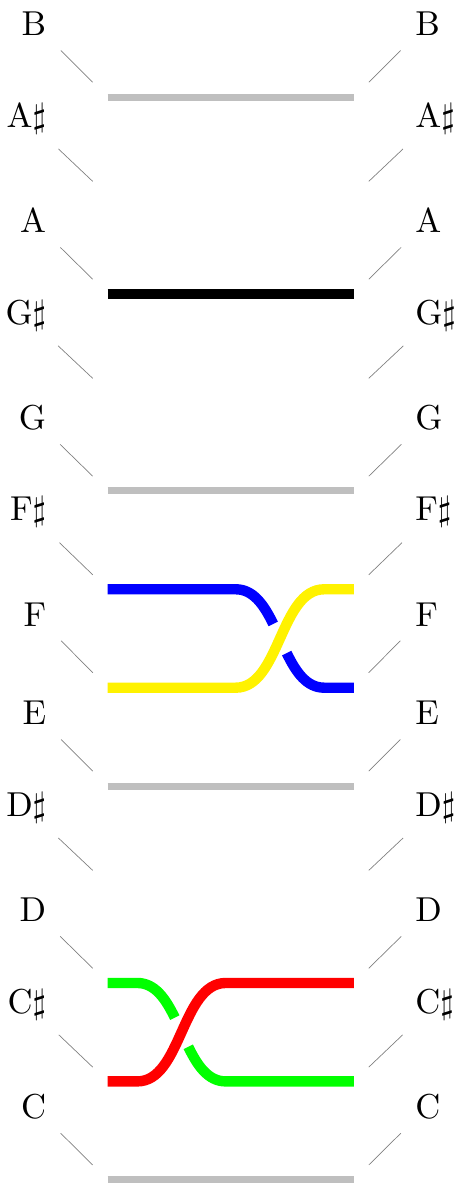}}
\hspace{15mm}
\subfloat[$C7$]
{\includegraphics[height=5cm]{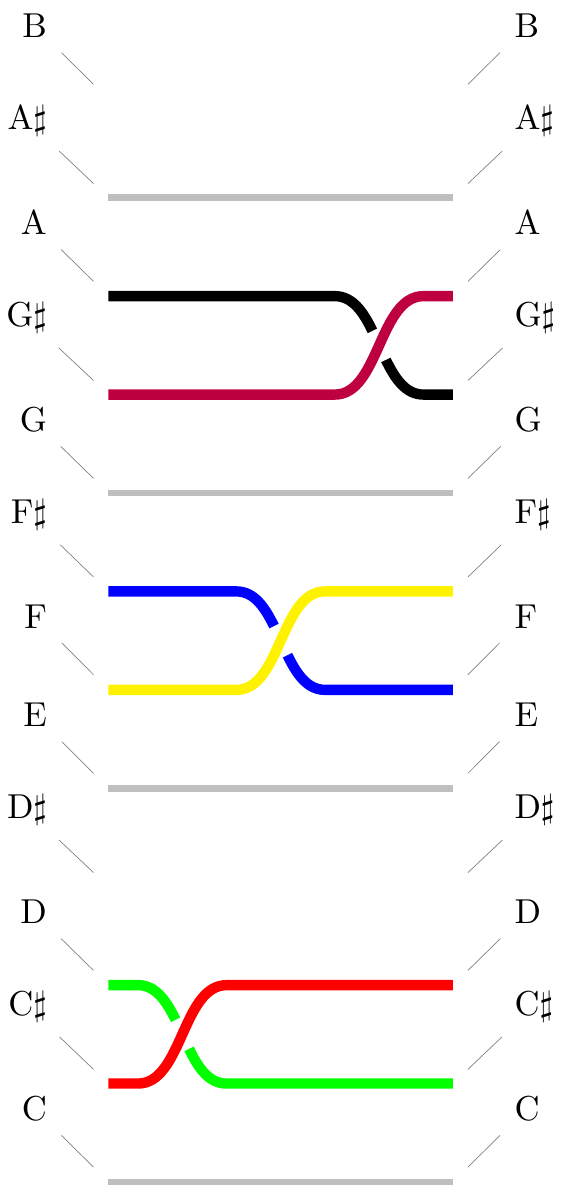}}
\hspace{15mm}
\subfloat[$C-7$]
{\includegraphics[height=5cm]{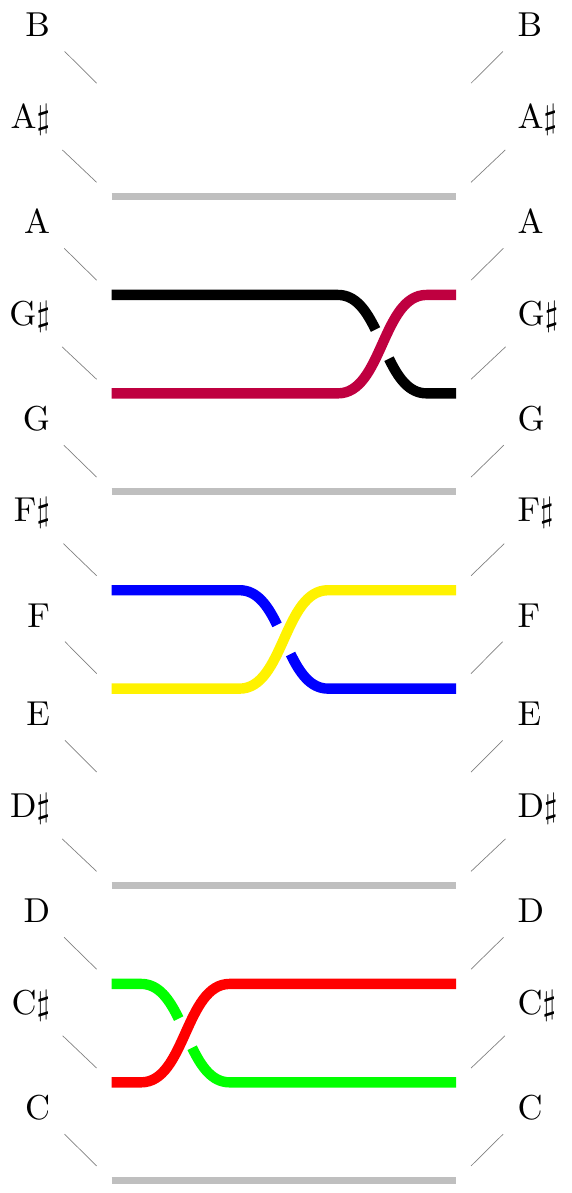}} 
\hspace{15mm}
\subfloat[$C-7^{\flat 5}$]
{\includegraphics[height=5cm]{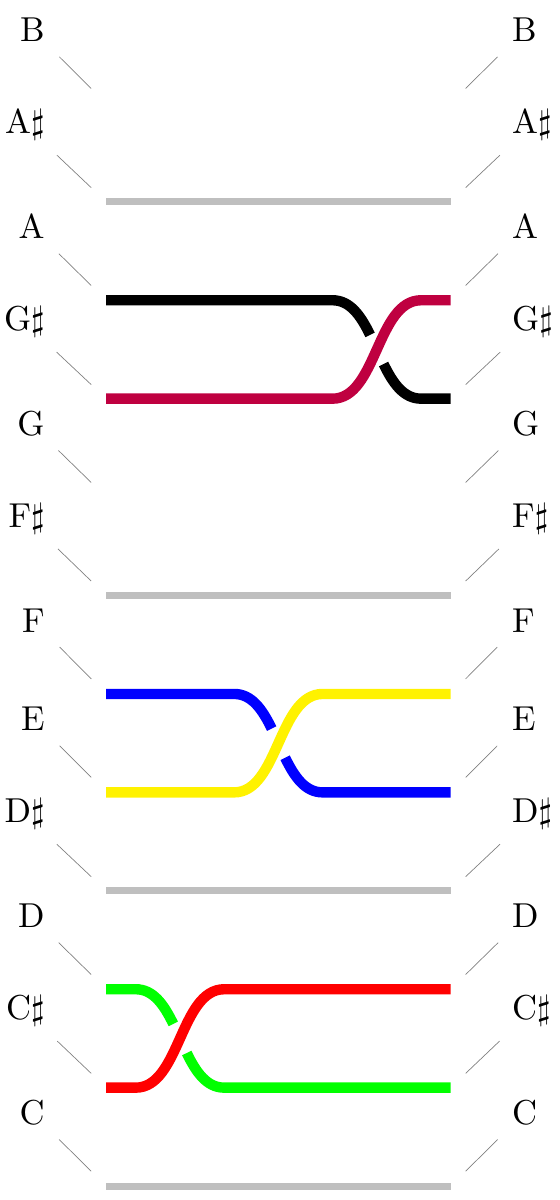}}
\caption{From graph to braids. These braids can be associated to the graphs we
represented in section \ref{sec:grafi}. Every loop in the graphs
turns out to be the represented by $\sigma_i$, see figure \ref{fig:braidgen}.}
\end{figure}

\subsection{A representation of cadential chord progression through
braids}\label{subsec:cad}

By using the concatenation product in the braid group, we shall represent 
harmonic progressions \cite[Chapter 12]{Pis59} through braids.

Given two base-chords $[B_1]$ and $[B_2]$, we can treat them as multisets of pitches, we denote the voice-leading among them
as $[B_1]\rightarrow [B_2]$ and we assume the voice leading to be crossing free. We refer to 
\cite{CQT08,Tym06, Tym09, Tym11, Pis59, Pis47, Pis55}.
Set $[B_1]=[(p_1,p_2,p_3,p_4)]$ and $[B_2]=[(q_1,q_2,q_3,q_4)]$, then 
we can write every note of each base-chord as an element of $\mathbb{R}/12\mathbb{Z}$, and arrange the voices of each chord, such that
\begin{equation}\label{eq:voice-cross}
p_i>p_j\Rightarrow q_i\geq q_j \mbox{ for all }i,j\in\{1,\dots ,4\}
\end{equation}

\begin{figure}[H]
\begin{centering}
\subfloat[$\mathcal{B}(Cmaj7)$]
{\includegraphics[height=7cm]{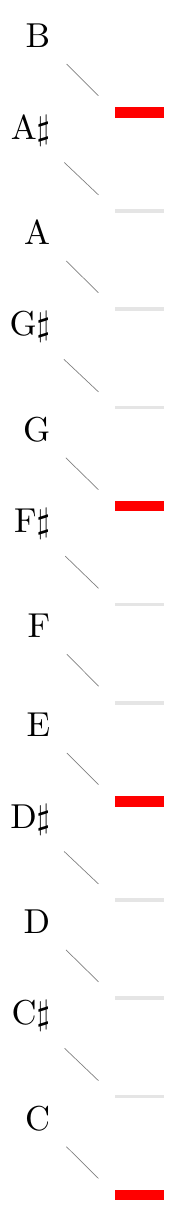}}
\hspace{15mm}
\subfloat[$\mathcal{B}(Gmaj7)$]
{\includegraphics[height=7cm]{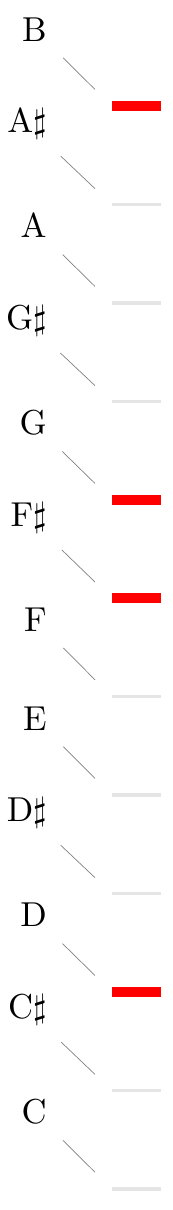}}
\hspace{15mm}
\subfloat[A voice moving on a braid]
{\includegraphics[height=7cm]{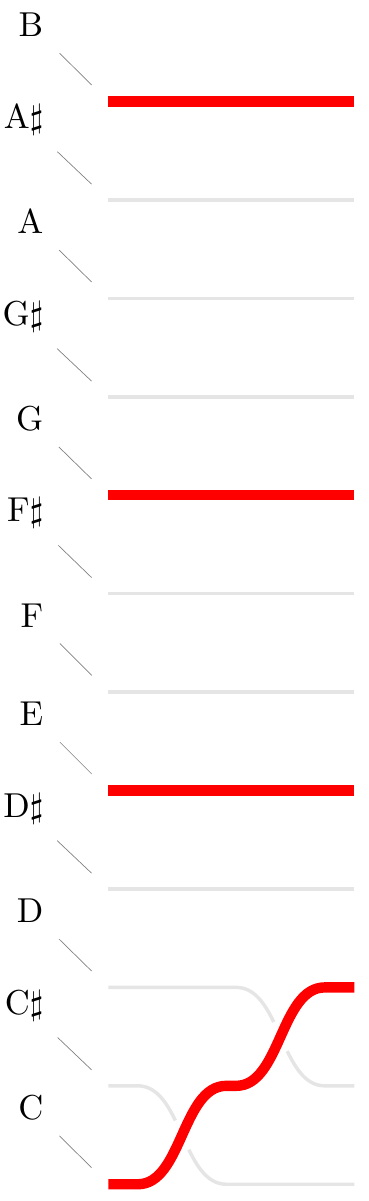}}
\caption{The braid diagram associated to the base-chords}
\label{fig:trivbraid}
\par\end{centering}
\end{figure}

\begin{example}
Consider $[Cmaj7]=[(C,E,G,B)]$ and $[Gmaj7]=[(G,B,D,F\sharp)]$ thus:
\[
[Cmaj7]=[(0,4,7,11)],\qquad
[Gmaj7]=[(7,11,2,6)].
\]

To avoid voice crossings in $(0,4,7,11)\rightarrow (7,11,2,6)$ we have
that
\[
 0 \mapsto 2, \ \ 4 \mapsto5, \ \ 7 \mapsto7, \ \ 11 \mapsto 11.
\]

which fulfills the condition (\ref{eq:voice-cross}).
\end{example}
A single chord can be represented as a trivial braid, we denote by $\mathcal B(\cdot)$ the braid associated to the chord $\cdot$. See diagrams (a) and (b) of figure
\ref{fig:trivbraid} for a braid representation of $[Cmaj7]$ and $[Gmaj7]$. 
The voice leading $Gmaj7\rightarrow Cmaj7$ generates a non trivial braid. For
example, to move from $[C]=[0]\in [Cmaj7]$ to $[D]=[2]\in [Gmaj7]$ (see figure
\ref{fig:trivbraid} (c)) we use the generators $\sigma_1\sigma_2$ (see section
\ref{sec:braids}). In this way it is possible to represent any 
chord progressions by braids. 

\subsection{Some jazz harmonic progressions}
In the following list we give a \emph{braid theoretic} representation of some 
of the most relevant jazz harmonic progressions.

\begin{itemize}
\item $II - V - I$ and $II - V - VI$.
The braid in figure \ref{fig:iiVi} (a), represents the harmonic progression
\[
D-7\rightarrow G7\rightarrow Cmaj7.
\]
where the second degree prepares the authentic cadence $V-I$.\\
In figure \ref{fig:iiVi} (b) is represented the harmonic progression
\[
D-7\rightarrow G7\rightarrow A-7
\]
where the second degree $D-7$ prepares the deceptive cadence $V-VI$
\begin{figure}[H]
\begin{centering}
\subfloat[$\mathcal{B}(D-7)\rightarrow \mathcal{B}(G7)\rightarrow
\mathcal{B}(Cmaj7)$]
{\includegraphics[width=5cm]{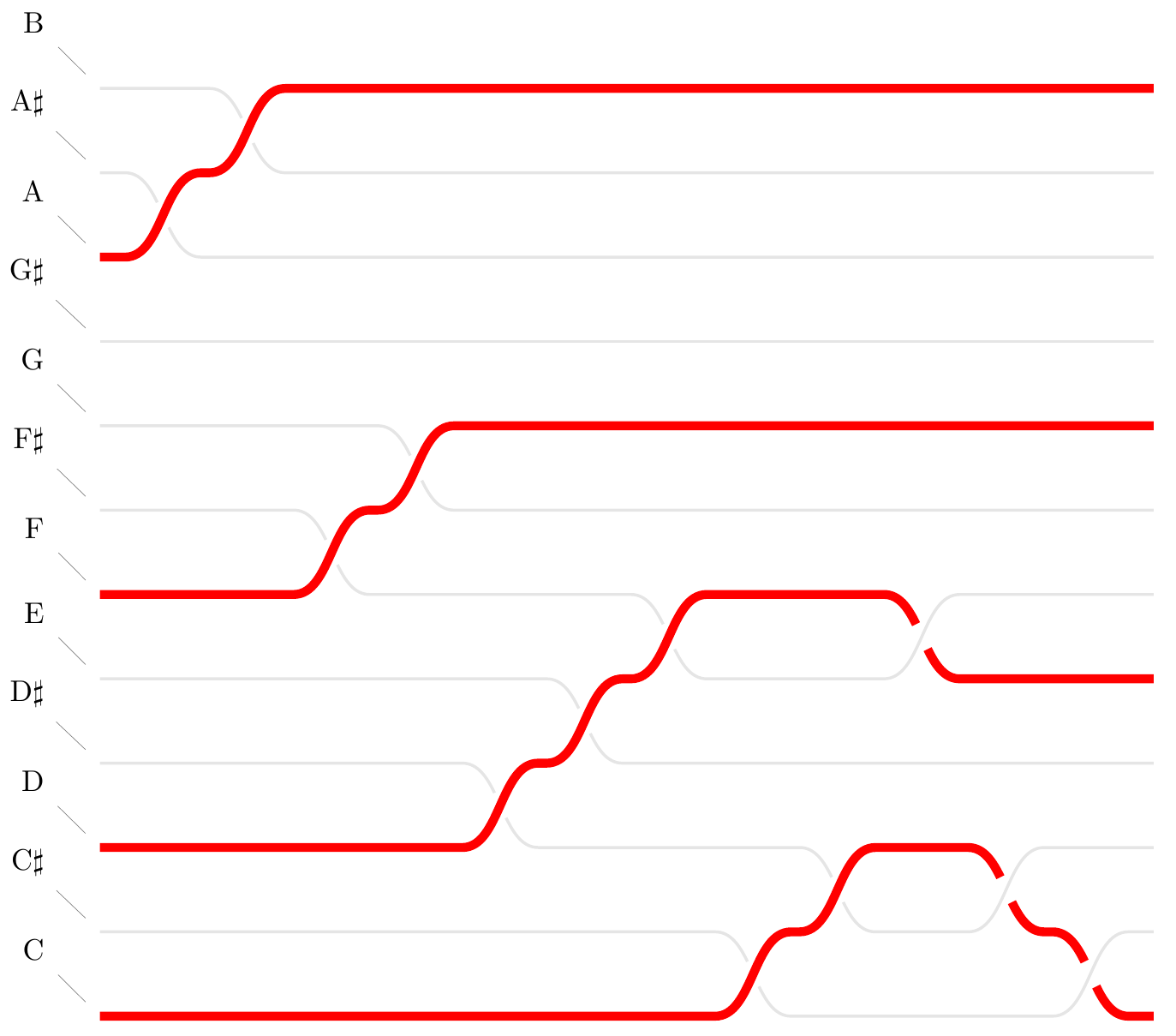}}
\hspace{15mm}
\subfloat[$\mathcal{B}(D-7)\rightarrow \mathcal{B}(G7)\rightarrow
\mathcal{B}(A-7)$]
{\includegraphics[width=6cm]{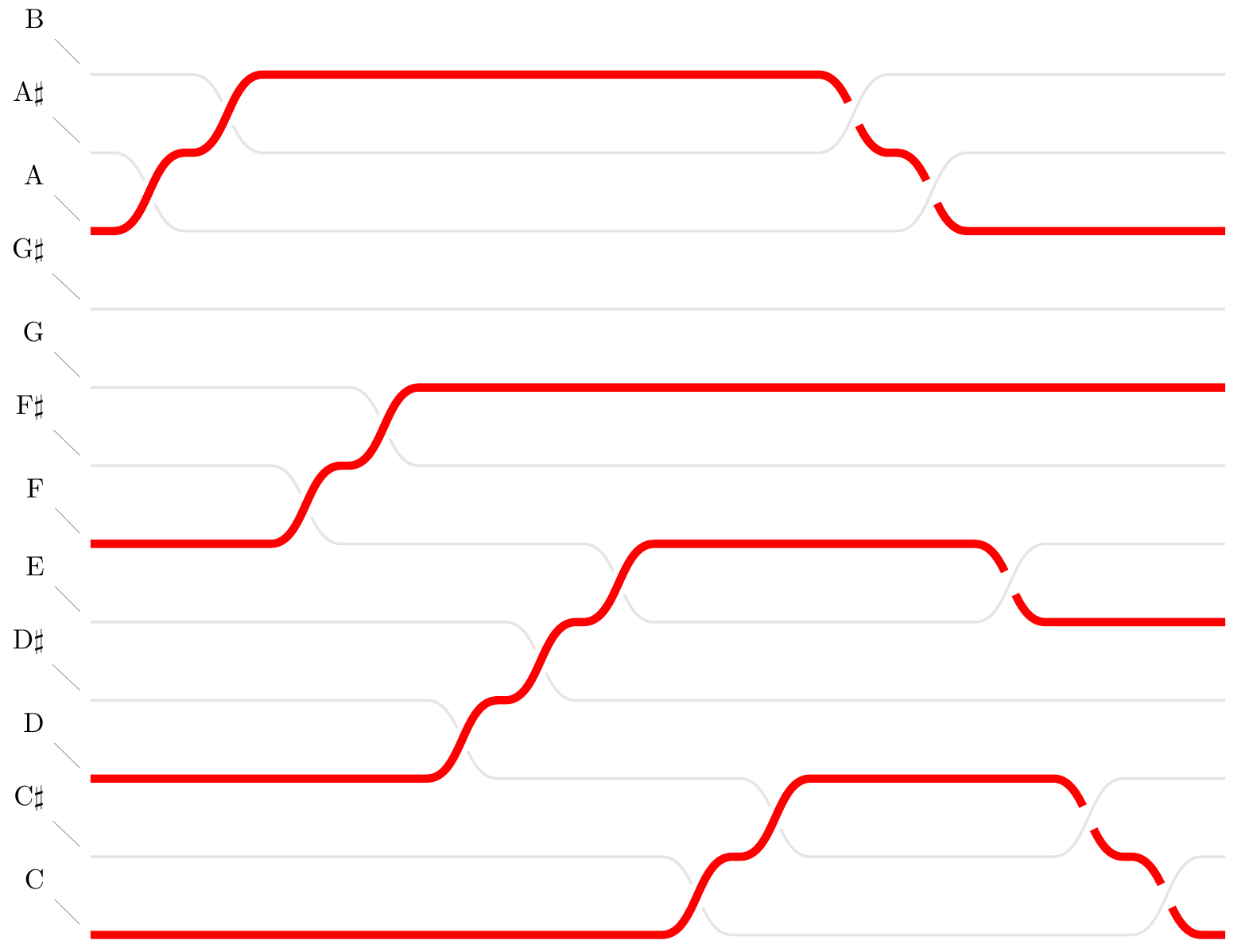}}
\caption{Authentic and deceptive cadences}
\label{fig:iiVi}
\par\end{centering}
\end{figure}

\item $IV - I$ and \emph{Secondary Dominant}. The braid in figure \ref{fig:IV-I}
(a) represents the plagal cadence
\[
Fmaj7\rightarrow Cmaj7.
\]
In figure \ref{fig:IV-I} (b) it is possible to see the harmonic progression
\[
F7\rightarrow C7\rightarrow G7
\]

\begin{figure}[H]
\begin{centering}
\subfloat[$\mathcal{B}(Fmaj7)\rightarrow \mathcal{B}(Cmaj7)$]
{\includegraphics[height=6.2cm, width=4cm]{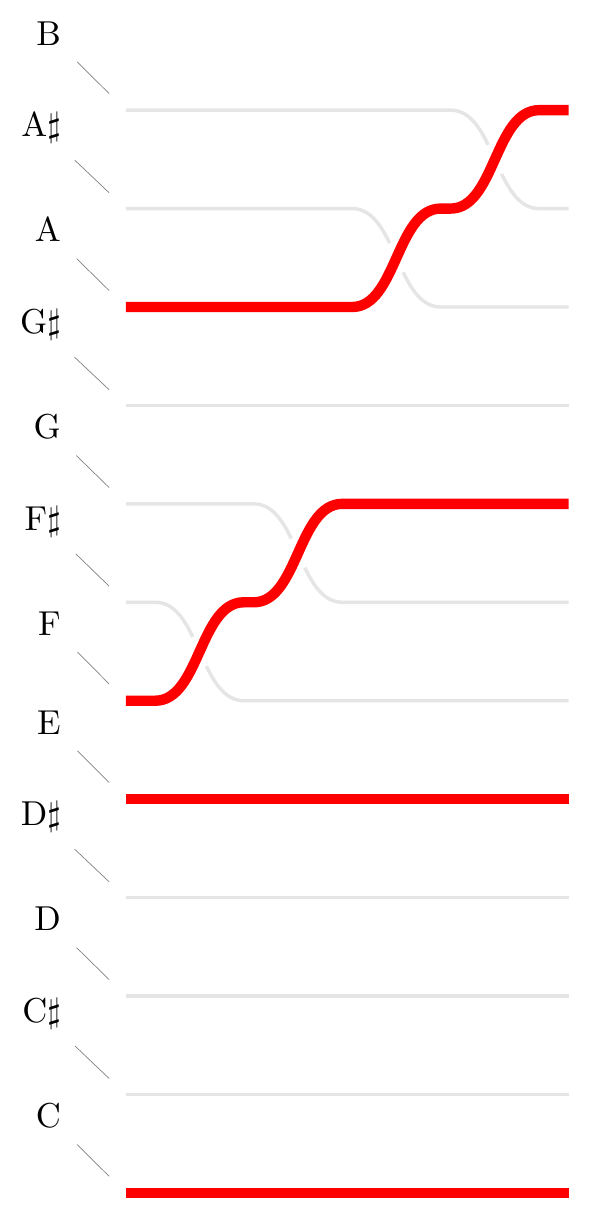}}
\hspace{15mm}
\subfloat[$\mathcal{B}(F7)\rightarrow \mathcal{B}(C7)\rightarrow
\mathcal{B}(G7)$]
{\includegraphics[width=5cm]{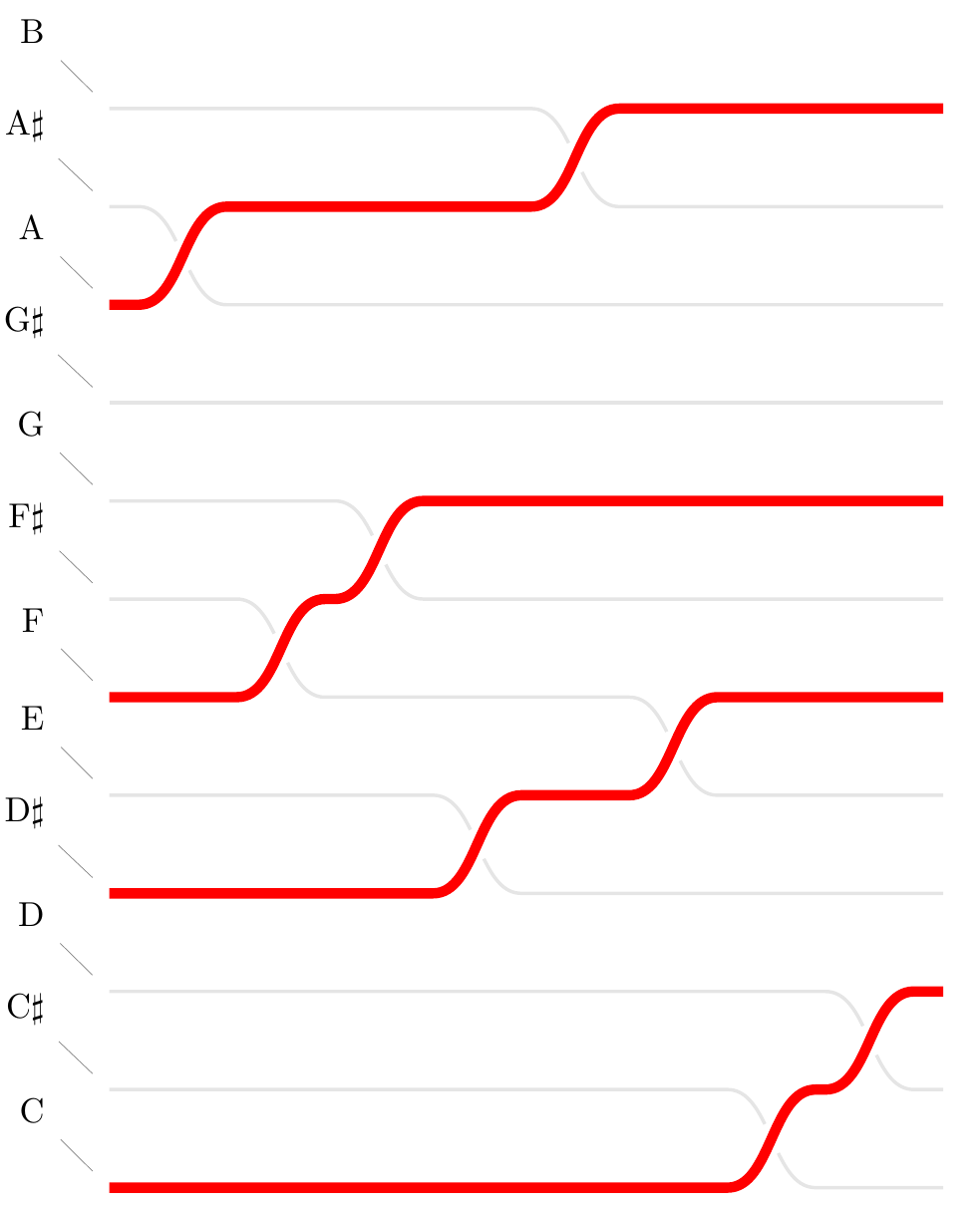}}
\caption{Plagal cadence and secondary dominant}
\label{fig:IV-I}
\par\end{centering}
\end{figure}
\end{itemize}


\section{An application: Mixolydian $\flat 2\,\sharp 4$ in Peru}
\begin{figure}[H]
\begin{centering}
\includegraphics[width=12cm]{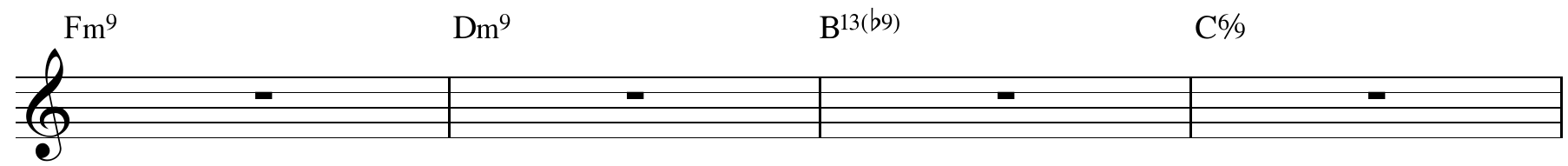}
\caption{An harmonic fragment from \emph{Peru }by Tribal Tech.}
\label{fig:harm}
\par\end{centering}
\end{figure}

\begin{figure}[H]
\centering
\subfloat[$\mathcal{B}(F-9)\rightarrow \mathcal{B}(D-9)$]
{\includegraphics[width=5cm, height=5cm]{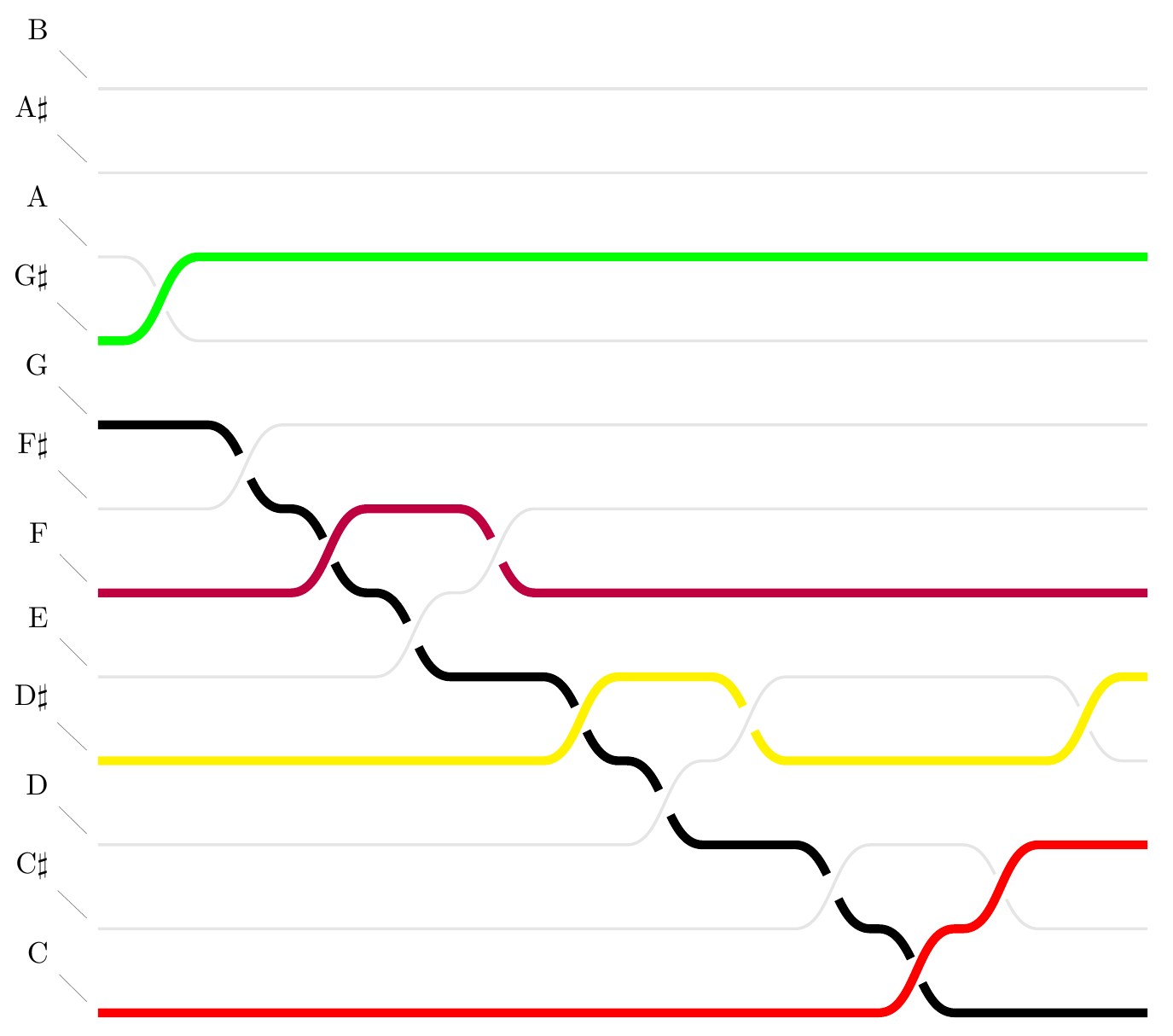}}
\hspace{10mm}
\subfloat[$\mathcal{B}(D-9)\rightarrow \mathcal{B}(B13^{\flat 9})$]
{\includegraphics[width=5cm, height=5cm]{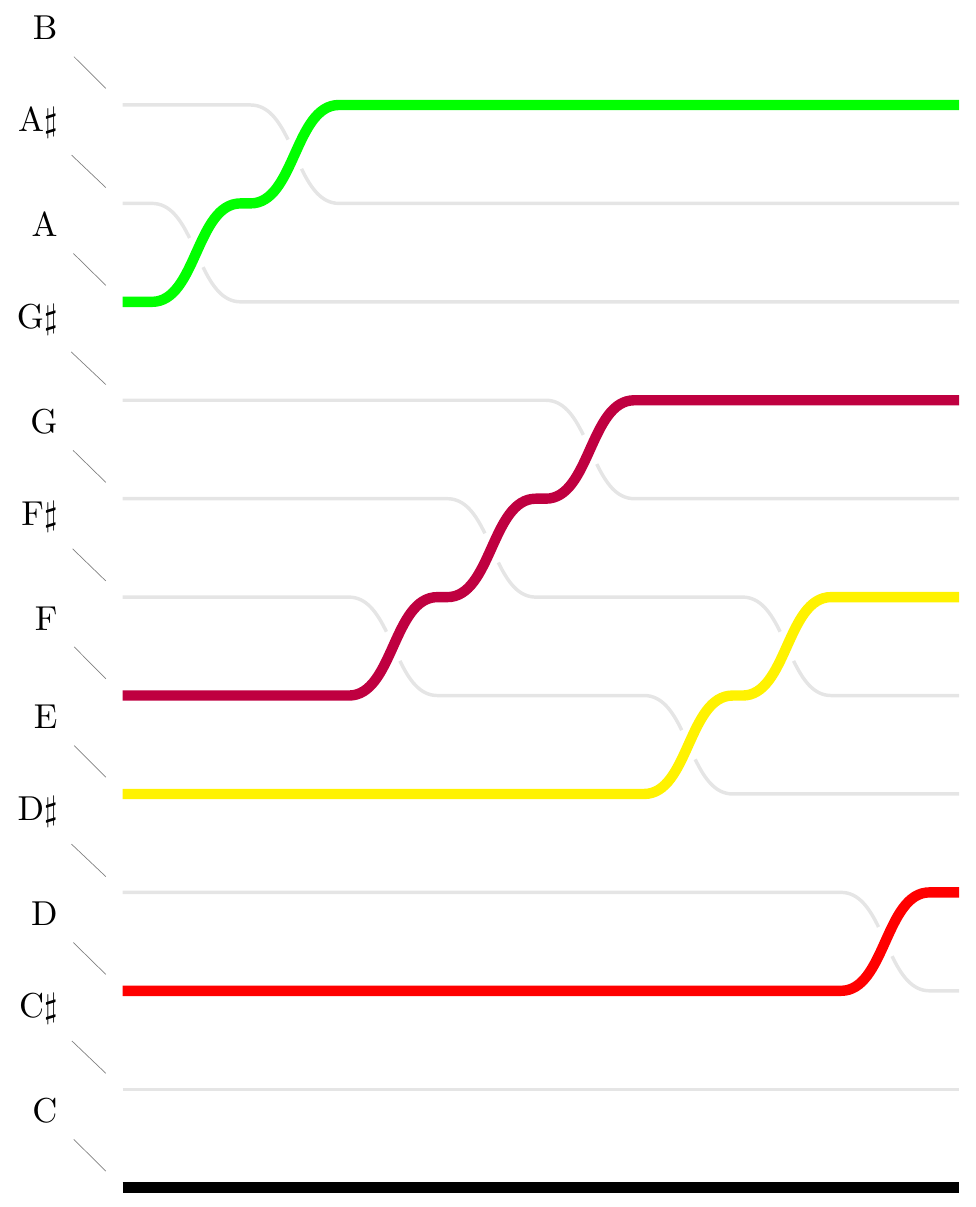}}
\vspace{5mm}
\subfloat[$\mathcal{B}(B13^{\flat 9})\rightarrow \mathcal{B}(C6/9)$]
{\includegraphics[width=5cm, height=5cm]{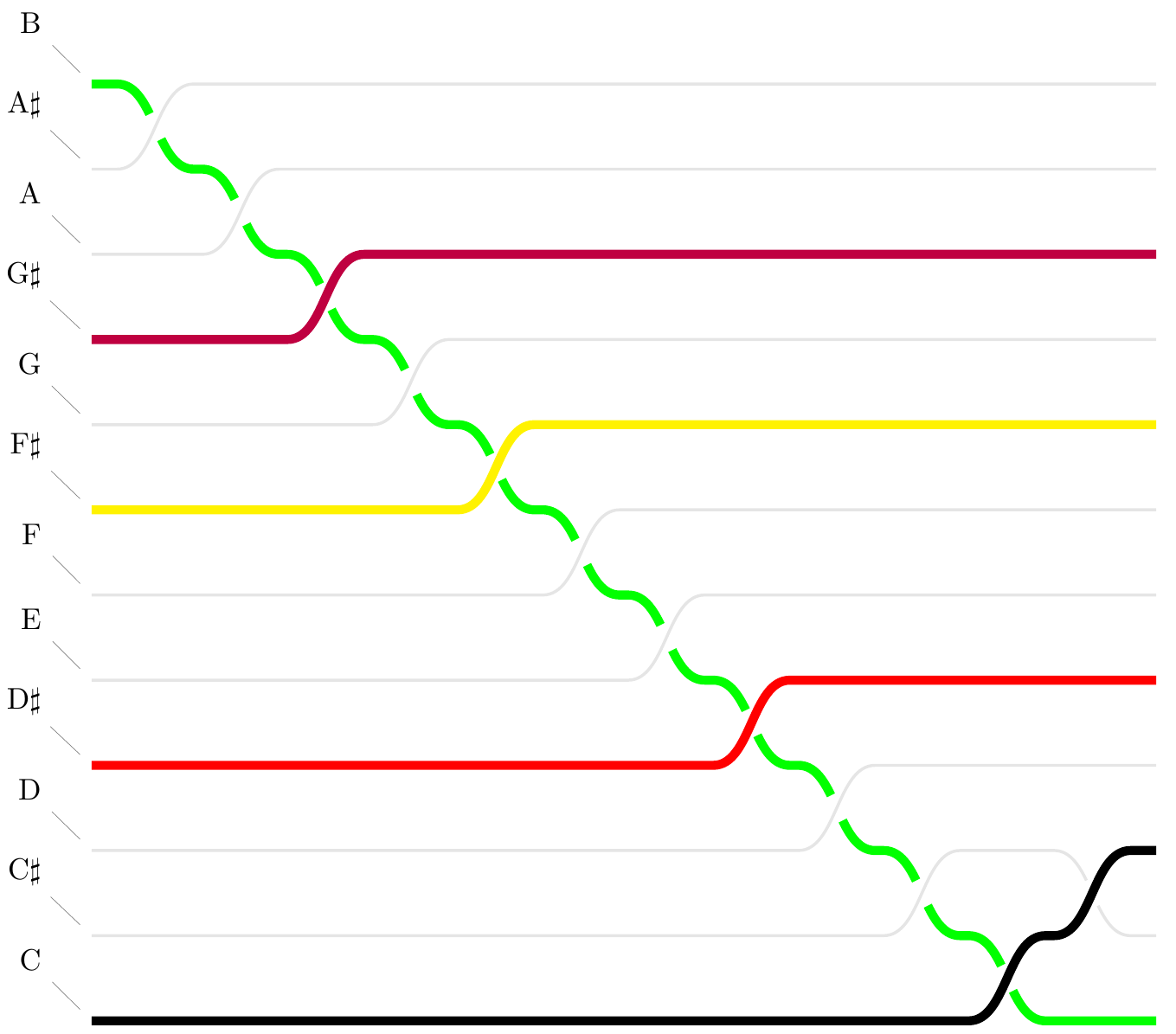}}
\caption{Voice leadings}
\label{fig:scott1}
\end{figure}

We introduced in paragraph \ref{subsec:cad} the representation through braids of chords
progressions. Let us analyse a modal harmonic structure excerpted from
\emph{Peru} by \emph{Tribal Tech} (see figure \ref{fig:harm}).
This choice is due to the clear harmonic and modal analysis of this track which
is provided by Scott Henderson in \cite{Hen}.
First of all we give a representation through braids of the harmonic progression. The algorithm we used to draw the braids which represent
voice-leadings between chords is exactly the one we described in
\ref{subsec:cad}. See figure \ref{fig:scott1} for such representation.

The braid \ref{fig:scott1} (a) is generated by
\[
 \sigma_9 \sigma_7 \sigma_6 \sigma_5 \sigma_6 \sigma_4 \sigma_3 \sigma_4
\sigma_2 \sigma_1 \sigma_2 \sigma_4
\] 

While, braids in figure \ref{fig:scott1} (b) and (c) have generators 
\[
\sigma_{10} \sigma_{11} \sigma_6 \sigma_7 \sigma_5 \sigma_6 \sigma_3
\]
and 
\[
 \sigma_{11} \sigma_{10}  \sigma_9 \sigma_8 \sigma_7 \sigma_6 \sigma_5 \sigma_4
\sigma_3 \sigma_2 \sigma_1 \sigma_2
\]
respectively.
Either an improviser or a composer can see each chord of the structure in figure
\ref{fig:harm} as a base-chord to build a modal melody. In \cite{Hen} the
following modal choices has been used to build a solo:
\begin{table}[H]
\begin{centering}
\begin{tabular}{l l l l}
$F-9$ & $D-9$ & $B^{13\flat 9}$ & $C6/9$ \\
\hline
F dorian & D dorian & B octatonic (half step-whole step) & C
lydian\tabularnewline
\end{tabular}
\par
\end{centering}
\begin{centering}
\caption{Example of modal choices on a given harmonic structure}
\label{tabperu}
\par\end{centering}
\end{table}
We do not considered any octatonic scale, however, thanks to the special modes
we introduced in proposition \ref{prop:specmodes}, we can give an heptatonic
approximation of the \emph{half step, whole step octatonic scale} (hs-ws scale
for brevity). See table \ref{tab:octa} for the notes of the octatonic scale with
root note $B$. 
\begin{table}[H]
\begin{centering}
\begin{tabular}{cccccccc}
$\mathbf{B}$ & $C$ & $D$ & $\mathbf{D\sharp}$ & $F$ & $\mathbf{F\sharp}$ &
$G\sharp$ & $\mathbf{A}$ \tabularnewline
\end{tabular}
\par
\end{centering}
\begin{centering}
\caption{Octatonic scale. Bold notes belong to a $B7$ chord}
\label{tab:octa}
\par\end{centering}
\end{table}

Since the base-chord for the hs-ws scale is a dominant seventh chord, it is
reasonable to look for a good approximation of this scale among the modes which
are built on a $B7$:
\begin{itemize}
\item B mixolydian: $(B,C\sharp,D\sharp,E,F\sharp,G\sharp,A)$
\item B lydian dominant: $(B,C\sharp,D\sharp,E\sharp,F\sharp,G\sharp,A)$
\item B mixolydian $\flat 13$: $(B,C\sharp,D\sharp,E,F\sharp,G,A)$
\item B mixolydian $\flat 2\flat 13$: $(B,C,D\sharp,E,F\sharp,G,A)$
\item \emph{B mixolydian} $\flat 2$: $(B,C,D\sharp,E,F\sharp,G\sharp,A)$
\item \emph{B mixolydian} $\flat 2\sharp 4$:
$(B,C,D\sharp,E\sharp,F\sharp,G\sharp,A)$
\item \emph{B mixolydian} $\sharp 4\flat 13$:
$(B,C\sharp,D\sharp,E\sharp,F\sharp,G,A)$
\item \emph{B mixolydian} $\flat 2\sharp 4 \flat 13$:
$(B,C,D\sharp,E\sharp,F\sharp,G,A)$
\end{itemize}
where special modes are emphasized. The choice is almost trivial, to find the
best approximation of the hs-ws scale it suffices to consider the nearest
special mode to the chord $B^{13\flat 9}$, which shares as many notes as possible with the hs-ws scale, i.e. the mixolydian $\flat 2\sharp 4$.
The comparison between the two scale is shown in table \ref{tab:comparison}, where the black square represents the note missing in the 
approximation of the 8-note scale with a 7-note one. The dropped note respect to the root of the base-chord is the minor third, which is a 
huge tension for a dominant chord, however the approximation provided by the special mode is a good one, since we are not giving up the tensions 
that are explicit in the chord $B^{13\flat 9}$ and at the same time the modal scale contains the arpeggio of the $B7$ chord.
\begin{table}[H]
\begin{centering}
\begin{tabular}{ccccccccc}
hs-ws scale & $B$ & $C$ & $D$ & $D\sharp$ & $F$ & $F\sharp$ &
$G\sharp$ & $A$ \tabularnewline
$mix\flat 2\sharp 4$ scale & $B$ & $C$ & {\large $\blacksquare$} & $D\sharp$ & $E\sharp = F$ & $F\sharp$ & $G\sharp$ & $A$ \tabularnewline
\end{tabular}
\par
\end{centering}
\begin{centering}
\caption{Comparing the hs-ws scale and the most suitable special mode.}
\label{tab:comparison}
\par\end{centering}
\end{table}


\section{Conclusion and future projects}\label{sec:conclusions}

We start this section by putting on evidence the major achievements of this paper. We  
starting by revising an heptatonic modal scale as a superimposition of two chords, playing 
two completely different roles. A base-chord $[B]$ which is represented by a seventh chord 
and a tension-triad constructed above it. Thus we can summarise it by as: ``$\,7=4+3.$'' From a  
genuinely musical viewpoint this introduce a striking difference since magnify the 
harmonic and melodic characteristic of the modal scale. 

The second achievement is to associate a 2-dimensional oriented planar graph to each base-chord 
reflecting its freedom to support tensions. By using this graph we are able to pull out from 
the same tonal system, or even better on the classical harminization, 
some {\em special modes\/} and hence some {\em new\/} (heptatonic) modal scales. 

The third, from our point of view, interesting achievement is the topological measure of 
complexity. This is given in terms of the Euler characteristic of a graph, and hence it is 
stable under homeomorphism of graphs. This is an important property making our construction 
independent on the choices of the positions of the degrees of the scale. From a mathematical point of 
view the Euler characteristic is an integer between $-2$ and $1$ included. Bigger this invariant, less 
degrees of freedom a base-chord has to support tension.  From a musical point of view this reflects a 
major degree of freedom in the improvisation as well as composing. 

Although this qualitative analysis give us a new insight on the modal scales, it is still poor in 
order to classify some modern tracks as well as to analyse some popular modern chord progressions and 
cadences. For this reason in section \ref{sec:braids} we introduce a braids theoretical interpretation 
of the above introduced graphs. This in particular allows us to associate to each chord progression, a 
braid in a quite effective and explicitly way by showing the voice-leading of the progression. 

We finally applied our results in order to better understand  \emph{Peru}. More precisely, by using the graph theoretic results 
and in particular the new special mode $\flat 2\,\sharp 4$ we are able to give a ``nice'' approximation of the octatonic scale by using an 
heptatonic one. Of course we have to give up to one note, but as we showed in the previous section, the 
note dropped out does not contribute to the tensions required explicitly by the chord. In fact the modal scale used as approximation 
includes the tensions which characterize the octatonic scale and the stable notes needed to be as melodic as possible on a such altered chord,
making our result very reasonable also from the musical point of view. We point out this approximation is 
good with respect to any chosen metric (intervallic, etc), simply because it shares 7 notes with the octatonic scale chosen by Scott Henderson. 

Using the braid theoretical approach to analyse Peru, we are able to figure out how complicate is 
the chord progression used by Scott Henderson in comparison to the most standard and popular 
chord progressions. 

In a forthcoming paper we shall try to distinguish about different melodic lines, constructed over the same base-chord progressions. 
From a mathematical viewpoint we shall construct some different topological spaces reflecting 
the crossings between the strand associated to a melodic line with respect to the underground braid. This 
will give us a sort of measure of complexity associated to a melodic phrase or to an improvisation and 
even to a classification of jazz and modern standard tracks.


\vspace{1cm}
Alessandro Portaluri\\
Department of Agriculture, Forest and Food Sciences\\
Universit\`a degli studi di Torino\\
Via Leonardo da Vinci, 44\\
10095 Grugliasco (TO)\\
Italy\\
E-mail: alessandro.portaluri@unito.it

\vspace{1cm}
\textcolor{white}{a}\\
Mattia G. Bergomi\\
Department of Computer Science\\
Universit\`a degli studi di Milano\\
Via Comelico, 39\\
20135 Milano (MI)\\
Italy\\
E-mail: mattia.bergomi@unimi.it

\end{document}